\magnification=\magstep1
\input amstex
\documentstyle{amsppt}
\catcode`\@=11 \loadmathfont{rsfs}
\def\mycal{\mathfont@\rsfs}
\csname rsfs \endcsname \catcode`\@=\active

\vsize=6.5in

\topmatter 
\title On the vanishing cohomology  problem \\ for cocycle actions of groups on II$_1$ factors  
\endtitle
\author SORIN POPA \endauthor

\rightheadtext{Vanishing cohomology}

\affil     {\it  University of California, Los Angeles} \endaffil

\address Math.Dept., UCLA, Los Angeles, CA 90095-1555, popa\@math.ucla.edu \endaddress

\thanks Supported in part by NSF Grant DMS-1700344 \endthanks

\abstract  We prove that any free cocycle action of a countable amenable group $\Gamma$ 
on  any II$_1$ factor $N$  
can be perturbed by inner automorphisms to a genuine action. Besides being satisfied by all amenable groups, 
this universal {\it vanishing cohomology} 
property, that we call $\Cal V\Cal C$, 
is also closed to free products with amalgamation over finite groups. While no other examples of $\Cal V\Cal C$-groups are known, 
by considering special cocycle actions $\Gamma \curvearrowright N$ in the case $N$ is the hyperfinite II$_1$ factor $R$,  
respectively the free group factor $N=L(\Bbb F_\infty)$, 
we  exclude many groups from being $\Cal V\Cal C$. 
We also show that any free action $\Gamma \curvearrowright R$ gives rise to a free cocycle $\Gamma$-action 
on the II$_1$ factor $R'\cap R^\omega$ whose vanishing cohomology is equivalent to Connes' Approximate Embedding 
property for the II$_1$ factor $R\rtimes \Gamma$. 
\endabstract

\endtopmatter

\document

\heading 0. Introduction \endheading

A {\it cocycle action} of a group $\Gamma$ on a II$_1$ factor $N$ is a map $\sigma: \Gamma \rightarrow \text{\rm Aut}(N)$ 
which is multiplicative modulo inner automorphisms of $N$, 
$\sigma_g \sigma_h= \text{\rm Ad}(v_{g,h}) \sigma_{gh}$, $\forall g,h\in \Gamma$, with the unitary elements 
$v_{g,h}\in \Cal U(N)$ satisfying the cocycle relation $v_{g,h} v_{gh,k} = \sigma_g (v_{h,k}) v_{g,hk}$, $\forall g,h,k\in \Gamma$.

If $\Gamma$ is a free group $\Bbb F_n$, for some $1\leq n \leq \infty$, 
then any cocycle $\Gamma$-action on any II$_1$ factor $N$ can obviously be perturbed by inner automorphisms $\{\text{\rm Ad}(w_g)\}_g$  
of $N$ so that to become a ``genuine'' action, i.e., such that $\sigma'_g = \text{\rm Ad}(w_g)\sigma_g$ 
is a group morphism, in fact so that the stronger condition $v_{g,h}=\sigma_g(w_h^*)w_g^*w_{gh}$, $\forall g,h$, holds true.   
We obtain in this paper several results towards identifying 
the class $\Cal V\Cal C$ of all countable groups $\Gamma$ that satisfy this 
universal vanishing cohomology property. Thus, we first 
prove that any free product of amenable groups amalgamated over a common finite subgroup 
is in the class $\Cal V\Cal C$. Then we show that 
if a group $\Gamma$ has an infinite subgroup which either has relative property (T), 
or has non-amenable centralizer, then $\Gamma$ is not in $\Cal V\Cal C$. To prove that all amenable groups lie in $\Cal V\Cal C$   
we use subfactor techniques to reduce the problem to the case $N$ is the hyperfinite II$_1$ factor $R$, where 
vanishing cohomology holds due to results in ([Oc85]). To exclude groups from being in $\Cal V\Cal C$ we apply 
W$^*$-rigidity results  to two types of cocycle actions that are ``hard to untwist'':   
the ones arising from $t$-amplifications of Bernoulli actions on $N=R$ introduced in [P01];  
and the ones considered in [CJ84], arising from normal inclusions $\Bbb F_\infty \hookrightarrow \Bbb F_n$ 
with $\Bbb F_n/\Bbb F_\infty = \Gamma$, which give  
cocycle $\Gamma$-actions on $N=L(\Bbb F_\infty)$.  

Untwisting cocycle actions on II$_1$ factors is a basic question 
in non-commutative ergodic theory and very specific to this area.  
Besides its intrinsic interest, the problem occurs in the classification of group actions 
on II$_1$ factors ([C74], [J80], [Oc85], [P01a]) and, closely related to it, in the classification of factors through unique crossed-product decomposition 
(as in [C74] for amenable factors, or [P01a], [P03], [P06a], [IPeP05], [PV12] for non-amenable II$_1$ factors). 
Another aspect, which goes back to ([CJ84]) and is important in W$^*$-rigidity, relates non-vanishing cohomology  
for certain cocycle $\Gamma$-actions on $L(\Bbb F_\infty)$ to non-embeddability of $L(\Gamma)$ into $L(\Bbb F_n)$. 
From an opposite angle, which offers a new point of view much emphasized here,   
vanishing cohomology results for cocycle actions are relevant to embedding problems, such as 
finding unusual group factors that embed into $L(\Bbb F_2)$ and Connes Approximate Embedding conjecture.

To describe the results in this paper in more details we need some background and notations. 
Let us first note that cocycle actions are more restrictive than {\it outer actions}, 
which are maps $\sigma: \Gamma \rightarrow \text{\rm Aut}(N)$  
that only require $\sigma_g\sigma_h \sigma_{gh}^{-1}\in \text{\rm Inn}(N)$, $\forall g,h\in \Gamma$. 
It has in fact been shown in ([NT59]) that there is a scalar 3-cocycle $\nu_\sigma  
\in H^3(\Gamma)$ associated to an outer action $\sigma$. If $\sigma$ is {\it free}, 
i.e.,  $\sigma_g\not\in \text{\rm Inn}(N)$, $\forall g\neq e$, then $\nu_\sigma$ is trivial if and only if $\sigma$ is a cocycle action. 
Thus, if we view the vanishing cohomology problem as 
a question about lifting a 1 to 1 group morphism $\sigma: \Gamma \rightarrow \text{\rm Out}(N)$ 
to a group morphism into $\text{\rm Aut}(N)$, then the problem is not well posed unless one requires $\nu_\sigma\equiv 1$, 
i.e., that $\sigma$ defines a cocycle action.  

Like for genuine actions, one can associate 
to a cocycle action $\Gamma \curvearrowright^\sigma N$ a tracial crossed product von Neumann algebra $N\rtimes \Gamma$,  
with the freeness of $\sigma$ equivalent to the condition 
$N'\cap N\rtimes \Gamma=\Bbb C1$. Thus, if $\sigma$ is free then $N\subset M=N\rtimes \Gamma$ is an 
irreducible inclusion of  II$_1$ factors with the normalizer of $N$ in $M$ generating $M$ as a von Neumann algebra 
($N$ is {\it regular} in $M$). Conversely, any irreducible regular 
inclusion of II$_1$ factors $N\subset M$ arises this way, from a crossed product construction 
involving a free cocycle action (cf. [J80]). 

The crossed product framework allows an alternative formulation of vanishing cohomology. 
Thus, if $M=N\rtimes \Gamma$ denotes the crossed product II$_1$ factor 
associated with the free cocycle action $(\sigma, v)$ of $\Gamma$ on $N$, 
and we let $\{U_g\}_g\subset M$ denote the canonical unitaries 
implementing $\sigma$ on $N$, then the existence of $w_g\in \Cal U(N)$ such that $v_{g,h}=w_g\sigma_g(w_h)w_{gh}^*$, $\forall g,h$ 
(i.e., {\it vanishing cohomology} for $v$) 
amounts to $U'_g=w_gU_g$ being a $\Gamma$-representation. While the condition that 
$\sigma'_g=\text{\rm Ad}(U'_g)$, $g\in \Gamma$, is a genuine action (i.e., {\it weak vanishing cohomology} for $v$) 
amounts to the weaker condition that 
$\{U'_g\}_g$ is a projective $\Gamma$-representation. 

Given a II$_1$ factor $N$, 
we denote by $\Cal V\Cal C(N)$ (respectively $\Cal V\Cal C_w(N)$) 
the class of countable groups $\Gamma$ with the property that any free cocycle action 
of $\Gamma$ on $N$ satisfies the strong form (respectively weak form) of the vanishing cohomology. 
Also, we denote by $\Cal V\Cal C$ (respectively $\Cal V\Cal C_w$) the class of 
countable groups $\Gamma$ with the property that any free cocycle action 
of $\Gamma$ on any II$_1$ factor $N$ satisfies the strong form (respectively weak form) of the vanishing cohomology.

The class $\Cal V\Cal C$ contains all finite groups by ([J80], [Su80]) and all groups with polynomial growth by ([P89]). The first main result 
in this paper, which we prove in Section 2,  shows that in fact $\Cal V\Cal C$ contains all countable amenable groups. Since by [J80] all  
1-cocycles for actions of finite groups are co-boundaries, this allows to deduce that, more than just containing the free groups, all amalgamated free products 
of amenable groups over finite groups belong to $\Cal V\Cal C$.  

\proclaim{0.1. Theorem} The class $\Cal V\Cal C$ contains all countable amenable groups. Also,  
if $\{\Gamma_n\}_n$ is a sequence of groups in $\Cal V\Cal C$ and $K\subset \Gamma_n$ is a common finite subgroup, $n\geq 1$, 
then $\Gamma_1*_K \Gamma_2 *_K ...\in \Cal V\Cal C$. 
\endproclaim

To prove the first part of this result we show that any cocycle action $\sigma$ of a countable amenable group $\Gamma$ on a 
separable II$_1$ factor $N$ can be perturbed by inner automorphisms to a cocycle action $\sigma'$ 
that leaves invariant an irreducible hyperfinite subfactor $R\subset N$ with the additional property that  
$\sigma'_g\sigma'_h{\sigma'_{gh}}^{-1}$ 
are implemented by unitaries in $R$, $\forall g,h$, 
with $\sigma'$ still free when restricted to $R$ (see Theorem 2.1).  This reduces the vanishing cohomology problem to the case 
$N=R$, where one can apply the vanishing cohomology result in ([Oc85]) to finish the proof.  

To prove the existence of a ``large'' $R\subset N$ that's normalized by inner perturbations of $\sigma$ 
we use an idea introduced in ([P89]; cf. also 5.1.5 in [P91]),  of
translating the problem into the question of whether there exists    
a sub-inclusion of hyperfinite factors inside 
the ``diagonal subfactor'' $N\simeq N^\sigma\subset M^\sigma$ associated with $\sigma$, so that to have a non-degenerate commuting square 
satisfying a {\it strong smoothness} condition on higher relative commutants. This subfactor  
problem was solved in [P89] in the case $\Gamma$ is finitely generated with trivial Poisson boundary (e.g., with polynomial growth; see [KV83]), 
by constructing $R$ as a limit of relative commutants $P_n'\cap N$ of factors in a 
tunnel  $N \supset N_1... $, obtained by iterating the downward basic construction (in the spirit of [P91], [P93]).  

However, that construction depends crucially on the trivial Poisson boundary condition on $\Gamma$. We use 
here the amenability of $\Gamma$ alone to construct a more elaborate decreasing sequence of  subfactors $P_n\subset N$ with $P_n'\cap N \nearrow R$ 
``large'' in $N$, obtained through {\it reduction/induction} in Jones tunnels. 
In fact, this method allows us to obtain the existence of strongly smooth embedding of hyperfinite 
subfactors into any finite index subfactor $N\subset M$ with standard invariant $\Cal G_{N\subset M}$ {\it amenable} (in the sense of [P91], i.e., 
with its graph $\Gamma_{N\subset M}$ satisfying the Kesten-type condition 
$\|\Gamma_{N\subset M}\|^2=[M:N]$; see also [P93], [P94a], [P97] for other equivalent definitions). We in fact show that given any  
amenable C$^*$-category $\Cal G$ of endomorphisms  
of a II$_\infty$ factor $\Cal N$ (viewed here as an outer action of an abstract rigid C$^*$-tensor category), 
there exists a ``large'' approximately finite dimensional (AFD) 
II$_\infty$ subfactor $R\overline{\otimes} \Cal B(\ell^2 \Bbb N)$ in $\Cal N$  that's normalized by $\Cal G$, modulo inner automorphisms (see Theorem 2.12). 

In Section 5 we use the strong solidity of free group factors ([OP07]) to 
prove that in order for a group $\Gamma$ to satisfy the property that any of its actions on II$_1$ factors normalizes 
a hyperfinite subfactor (modulo inner automorphisms), $\Gamma$ must necessarily be amenable. The problem of whether this 
dichotomy still holds for subfactor standard invariants  and rigid C$^*$-tensor categories, remains open.

In turn, in Sections 3 and 4 we obtain a series of obstruction criteria for groups to belong to the classes   
$\Cal V\Cal C(R), \Cal V\Cal C(L(\Bbb F_\infty)), \Cal V\Cal C$,  summarized in the following: 

\proclaim{0.2. Theorem} $1^\circ$ If a countable group $\Gamma$ has an infinite subgroup which either has the 
relative property $\text{\rm (T)}$, or has non-amenable centralizer in $\Gamma$, then $\Gamma \not\in \Cal V\Cal C_w(R)$. 

$2^\circ$ Assume a countable group $\Gamma$ satisfies one of the following: $(a)$ $\Gamma$ does not have 
Haagerup property $($e.g., it contains an infinite subset with relative property $\text{\rm (T)})$; $(b)$ The Cowling-Haagerup 
invariant $\Lambda(\Gamma)$ is larger than $1$; $(c)$ $\Gamma$ has an infinite subgroup with non-amenable centralizer; $(d)$ $\Gamma$ has 
an infinite amenable subgroup with non-amenable normalizer. Then $\Gamma \not\in \Cal V\Cal C_w(L(\Bbb F_\infty))$,  
in particular $\Gamma \not\in \Cal V\Cal C_w$. 
\endproclaim 

To prove the restrictions on $\Cal V\Cal C(R)$ we use the $t$-amplifications of Bernoulli actions on $R$ introduced in [P01a] and 
results obtained there and in ([P06a]) through deformation-rigidity arguments. In turn, to get restrictions on  $\Cal V\Cal C_w(L(\Bbb F_\infty))$, 
we use the {\it Connes-Jones} (CJ) {\it cocycles} associated with surjective group morphisms $\pi: \Bbb F_S \rightarrow (\Gamma, S)$, extending  
the map assigning the free generators of $\Bbb F_S$ to a set of generators $S\subset \Gamma$.  As shown in [CJ84], 
if $\Gamma$ is infinite, non-free,  
then $\text{\rm ker}\pi\simeq \Bbb F_\infty$ and the inclusion $L(\Bbb F_\infty)=N \subset M = L(\Bbb F_S)$ is of the form 
$N\subset N\rtimes_{(\sigma_\pi,v_\pi)} \Gamma$, for a free cocycle action $(\sigma_\pi,v_\pi)$. The vanishing of the cocycle $v_\pi$ implies 
that $L(\Gamma)$ embeds into $L(\Bbb F_S)$, hence $2^\circ$ above follows from results in ([CJ84], [P01b], [O03], [P06b], [OP07]). 

The CJ-cocycles seem the ``most difficult to untwist'', in the sense that if all such cocycle actions of  
$\Gamma$ on $L(\Bbb F_\infty)$ untwist, then $\Gamma$ ought to be in $\Cal V\Cal C$. 
In particular, this would show that $\Cal V\Cal C = \Cal V\Cal C(L(\Bbb F_\infty))$. Since untwisting a CJ-cocycle for $\Gamma$ implies that 
$\Gamma$ is in the class 
W$^*_{leq}(\Bbb F_2)$ of groups whose von Neumann algebra embeds into $L(\Bbb F_2)$, one has 
$\Cal V\Cal C \subset \Cal V\Cal C(L(\Bbb F_\infty)) \subset \text{\rm W}^*_{leq}(\Bbb F_2)$. Very little is actually known about the class  $\text{\rm W}^*_{leq}(\Bbb F_2)$, 
which is extremely interesting by itself. Any success in proving $\Cal V\Cal C$ property for some ``exotic'' group $\Gamma$, would provide 
embeddings  $L(\Gamma) \hookrightarrow L(\Bbb F_2)$.   

In the final part of the paper we discuss a connection between vanishing cohomology phenomena and Connes Approximate Embedding  
conjecture, on whether any separable II$_1$ factor $M$ embeds in the ultrapower $R^\omega$ of the hyperfinite II$_1$ factor $R$. 
Thus, we notice that any free action of a group $\Gamma$ on $R$ (such as the ``non-commutative'' 
Bernoulli action $\Gamma \curvearrowright R^{\overline{\otimes}\Gamma} \simeq R$), 
gives rise to a free cocycle action of $\Gamma$ 
on the centralizer $R_\omega=R'\cap R^\omega$ of $R$ in $R^\omega$. 
We deduce that this cocycle untwists if and only $M=R\rtimes_\sigma \Gamma$ satisfies 
the  conjecture. 
 
{\it Acknowledgement.} I am very grateful to Damien Gaboriau, Vaughan Jones, Jesse Peterson and Stefaan Vaes for 
many useful discussions related to this paper. I am also grateful to the referee for his/her many 
pertinent questions that led to what I believe to be a much improved  final version. 

\heading 1. Preliminaries and notations. \endheading

For general background on II$_1$ factors we refer the reader to ([AP17]; also [T79], [BrO08] for general theory of 
operator algebras and von Neumann algebras). 
\vskip.05in

\noindent
{\bf 1.1. Cocycle actions and crossed products.} Given a II$_1$ factor $N$, we denote by 
Aut$(N)$ the group of automorphisms of $N$. An automorphism $\theta$ of $N$ is {\it inner} 
if there exists $u$ in the unitary group of $N$, $\Cal U(N)$, 
such that $\theta(x) = {\text{\rm Ad}}u (x)=uxu^*, \forall x\in N$. We denote by Inn$(N)\subset {\text{\rm Aut}}(N)$ 
the group of all such inner automorphisms and by $\text{\rm Out}(N)$ the quotient group $\text{\rm Aut}(N)/\text{\rm Inn}(N)$.  

Given a discrete group $\Gamma$, an {\it action} of $\Gamma$ on 
$N$ is a group 
morphism  $\sigma : 
\Gamma \rightarrow \text{\rm Aut}(N)$. We will use the notation $\Gamma\curvearrowright^\sigma N$ to emphasize such an action. 

More generally, a {\it cocycle action} $\sigma$ of $\Gamma$ on $N$ 
is a map $\sigma : \Gamma \rightarrow \ \text{\rm Aut} (N)$ 
with the property that there exists    
$v : \Gamma \times \Gamma \rightarrow \Cal U(N)$ such that:
$$
\align
& \sigma_e =\text{\rm id and} \ 
\sigma_g \sigma_h = \ \text{\rm Ad} \ v_{g,h} \sigma_{gh}, \ \forall 
g, h \in G \tag 1.1.1 \\
& v_{g,h} v_{gh,k} = \sigma_g (v_{h,k}) v_{g,hk}, \ \forall g, h, k \in \Gamma. 
\tag 1.1.2
\endalign
$$

The cocycle action $\sigma$ is {\it free} if  $\sigma_g$ 
cannot be implemented by unitary elements in $N, \forall g\neq e$, in other words if 
the factoring of $\sigma$ through the quotient map $\text{\rm Aut}(N) \rightarrow \text{\rm Out}(N)$ is 1 to 1. 
All cocycle actions  (in 
particular all actions) that we will consider in this paper are assumed to be free. 

Following ([KT02]), a map $\sigma: \Gamma \rightarrow \text{\rm Aut}(N)$ 
that's a 1 to 1 group morphism when factored through the quotient map $\text{\rm Aut}(N) \rightarrow \text{\rm Out}(N)$, 
is called an {\it outer $\Gamma$-action} (an alternative terminology for such $\sigma$ is $Q$-{\it kernel}, 
notably used in [J80], [Oc85]). 
Thus, an outer action satisfies $(1.1.1)$ above, but not necessarily $(1.1.2)$.  As shown in ([NT59]), 
if $\sigma$ is an outer $\Gamma$-action, then one can associate to it a scalar 3-cocycle $\nu_\sigma 
\in H^3(\Gamma)$, with the property that $\nu_\sigma\equiv 1$ if and only if $\sigma$ is a cocycle action, and which one calls 
the $H^3(\Gamma)$-{\it obstruction} of $\sigma$.

If $\sigma$ is a cocycle action, then a map $v$ satisfying (1.1.2) is called a {\it \text{\rm 2}-cocycle} for 
$\sigma$.  The 2-cocycle is {\it normalized} if $v_{g,e} = v_{e,g} = 1, 
\ \forall g \in G$.  Note that 
any 2-cocycle satisfies $v_{e,e} \in \Bbb C$. Thus any 2-cocycle 
$v$ can be normalized by replacing it, if necessary, by $v'_{g,h} = v^*_{e,e} \ 
v_{g,h}, \ g, h \in \Gamma$. All 2-cocycles considered from now on will be normalized.

Also, when given a cocycle action $\sigma$, we will 
sometimes specify from the beginning 
the 2-cocycle it comes with, thus considering it as a pair 
$(\sigma, v)$. 

Note that the $2$-cocycle $v$ is unique modulo perturbation by 
a {\it scalar} 2-{cocycle} $\mu$. More precisely,  
$v': \Gamma \times \Gamma \rightarrow \Cal U(N)$, 
with $v'_{e,e} = 1$, satisfies 
conditions (1.1.1), (1.1.2) if and only 
if $v'=\mu v$ for some  
scalar valued function $\mu : \Gamma \times \Gamma \rightarrow \Bbb T$ 
satisfying 
$\mu_{e,e} = 1$ and 
$$
\mu_{g,h} \mu_{gh,k} = \mu_{h,k} \mu_{g, hk}, 
\quad \forall g, h, k \in \Gamma 
\tag 1.1.3
$$

Let us recall the definition of the {\it crossed product} von Neumann  algebra associated 
with a cocycle action $(\sigma, v)$ of $\Gamma$ on $N$, denoted 
$N \rtimes_{(\sigma,v)} \Gamma$ (or simply $N\rtimes_\sigma \Gamma$ 
if $\sigma$ is a genuine action). So let $\Cal H$ denote the Hilbert space  
$\oplus_h (L^2N)_h\simeq \ell^2(\Gamma, \ L^2N)$, which we view as the space 
of $\ell^2$-summable formal series $\sum_h U_h\xi_h$, where  
$\{U_g\}_{g\in \Gamma}$ are here ``indeterminates'' (labels) and the ``coefficients'' $\xi_h$ belong to $L^2N$.

We define a $^*$-algebra structure on the subspace $\Cal H_0\subset \Cal H$ of finitely supported sums 
with ``bounded'' coefficients $\xi_g=x_g\in N$, and at the same time a Hilbert  
$\Cal H_0$-bimodule structure on $\Cal H$,  as follows. First, we let $U_e$ act on both left and right on $\Cal H$ as 
the identity $id_{\Cal H}$ and identify $N$ with $U_eN$ acting left-right on $\Cal H$ by $x (\sum_h U_h \xi_h)y=
\sum_h U_h (\sigma_h^{-1}(x)\xi_h y)$. Then we let  
$U_g\cdot \sum_h U_h \xi_h =\sum_h  v_{g,h}U_{gh} \xi_h$, which by the change of variables $h'=gh$ 
and ``moving'' $v_{g,h}$ from left to right according to the above multiplication by $N$ rule, is equal to  
$\sum_h U_{h} \sigma_{h}^{-1}(v_{g^{-1}h,h}) \xi_{g^{-1}h}$.  We also let $(\sum_h U_h \xi_h)\cdot U_g= 
\sum_h v_{h,g} U_{hg} \sigma_g^{-1}(\xi_h)$ which by similar considerations is equal to 
$\sum_h U_h \sigma^{-1}_h(v_{hg^{-1},g}) \sigma_g^{-1}(\xi_{hg^{-1}})$. Finally, we let $U_g^*=U_{g^{-1}}
v_{g,g^{-1}}^*$ and $(U_gx)^*= x^*U_g^*=U_{g^{-1}} v_{g,g^{-1}}^* \sigma_g(x^*)$. 
The $^*$-algebra $\Cal H_0$ has a trace given by $\tau(\sum_h U_h x_h)=\tau_N(x_e)$ which recovers 
for elements in $\Cal H_0$ the $\Cal H$-scalar product, i.e., if $X, Y\in \Cal H_0$ then 
$\langle X, Y \rangle_{\Cal H}=\tau(Y^*X)$.  

It is easy to verify that the left multiplication by $U_g$ give unitary operators on $\Cal H$, 
the left multiplication by $N=NU_e\subset \Cal H_0$ on $\Cal H$ gives a representation  
of $N$ as a von Neumann algebra,  and altogether left multiplication by elements in $\Cal H_0$ 
are bounded operators on $\Cal H$ that give a $^*$-representation $M_0$ 
of $\Cal H_0$ in $\Cal B(\Cal H)$, 
with the trace $\tau$ being implemented by the 
vector state given by $1=U_e1\in \Cal H$. 

The crossed product von Neumann algebra $N \rtimes_{(\sigma,v)} \Gamma$ is by definition 
the weak operator closure of $M_0$ in $\Cal B(\Cal H)$. It is a finite von Neumann algebra with faithful normal state 
$\tau(X)=\langle X1, 1 \rangle_{\Cal H}$, $\forall X\in M$. One clearly has a natural identification between the standard representation 
(or standard Hilbert $M$-bimodule) 
$L^2(M,\tau)$ and  $\Cal H$.

The cocycle action $(\sigma, v)$ is {\it free}  if $\sigma_g$ 
is an outer automorphism, $\forall g \neq e$. It is well known 
(and immediate to check) that $(\sigma, v)$ is free if and only if $N' \cap N 
\rtimes_{(\sigma,v)} \Gamma = \Bbb C1$. So in this case  
$M=N \rtimes_\sigma \Gamma$ is a II$_1$ factor with the normalizer $\Cal N_M(N)=\{u\in \Cal U(M) \mid uNu^*=N\}$ 
of $N$ in $M$ generating $M$ (i.e., with $N$ {\it regular} in $M$). 

Conversely, if $N\subset M$ is a regular, irreducible inclusion of II$_1$ factors and we denote $\Gamma = \Cal N_M(N)/\Cal U(N)$, 
with $U_g\in \Cal N_M(N), g\in \Gamma$, a lifting of $\Gamma$, $U_e=1$, and we denote $\sigma_g=\text{\rm Ad}(U_g)_{|N}$, 
$v_{g,h}=U_{gh}U_h^*U_g^*$, then $(\sigma, v)$ is a free cocycle action of $\Gamma$ on $N$, with $N\subset N\rtimes_{(\sigma,v)} \Gamma$ 
naturally isomorphic to $N\subset M$ (see e.g., [J80]). 

\vskip.05in

\noindent
{\bf 1.2.   Cocycle conjugacy of cocycle actions.} The cocycle actions 
$(\sigma_i, v_i)$ of 
$\Gamma$ on $N_i$, $i=1, 2$,  are {\it cocycle conjugate} if 
there exists an isomorphism 
$\theta:N_1 \simeq N_2$ and a map $w:\Gamma \rightarrow \Cal U(N_2)$ such that the following conditions are satisfied: 
$$
\theta \sigma_1(g) \ \theta^{-1} = \text{\rm Ad}(w_g) \sigma_2(g), \quad 
\forall g. \tag 1.2.1   
$$
$$
\theta(v_1(g,h)) = w_g \sigma_2(g) (w_h) v_2(g,h) \ w_{gh}^*, \ 
\forall g,h. \tag 1.2.2 
$$

The cocycle actions $\sigma_1, \sigma_2$  are {\it outer conjugate} (or 
{\it weakly cocycle conjugate}) if condition (1.2.1) 
is satisfied.  Note that this is equivalent to 
$\sigma_1, \sigma_2$ composed with the quotient map $\text{\rm Aut}(N)\rightarrow \text{\rm Out}(N)$  
being conjugate by an element in $\text{\rm Out}(N)$.

Similarly, two outer actions $\sigma_1, \sigma_2 \rightarrow \text{\rm Aut}(N)$ are {\it outer conjugate}, if 
there exists $\theta\in \text{\rm Aut}(N)$ such that $\sigma_1(g)=\theta \sigma_2(g)\theta^{-1}$ modulo $\text{\rm Inn}(N)$, 
for all $g\in \Gamma$. 

The actions $\sigma_1, \sigma_2$ are {\it conjugate} 
if there exists an isomorphism $\theta:N_1 \simeq N_2$ such that 
conditions (1.2.1) is satisfied 
with $w=1$. We then write $\sigma_1 \sim \sigma_2$. 

We recall here the following well known observation (see e.g., [J80]), 
which translates cocycle conjugacy of free cocycle actions into the isomorphism 
of the associated crossed-product inclusions of factors. 

\proclaim{Proposition} Let $(\sigma_i, v_i)$ be a cocycle action of 
the discrete group $\Gamma_i$ on the $\text{\rm II}_1$ factor $N_i$, $i=1,2$. If there exists a $*$-isomorphism 
$\Phi: N_1\rtimes_{(\sigma_1, v_1)} \Gamma_1 \simeq  N_2\rtimes_{(\sigma_2, v_2)} \Gamma_2$
such that $\Phi(N_1)=N_2$, then
$\sigma_1$ and $\sigma_2$ are cocycle conjugate. More precisely, 
there exists a group isomorphism $\gamma:\Gamma_1\rightarrow \Gamma_2$, and 
unitaries $w_g\in {\Cal U}(N_2)$, for all $g\in \Gamma_1$, such that:

(i) $\Phi \sigma_1(g) \Phi^{-1}={\text{\rm Ad}} w_g\text{ }\sigma_2({\gamma(g)})$, for 
all $g\in \Gamma_1$, 

(ii) $\Phi(v_1({g,h}))=w_g \sigma_2({\gamma(g)}(w_h)v_2({\gamma(g),
\gamma(h)})w^*_{gh}$, for all $g,h\in \Gamma_1$.

Conversely, if $\Phi:N_1\simeq N_2$ is a $*$-isomorphism, 
$\gamma:\Gamma_1\simeq \Gamma_2$ is a group isomorphism, and there exist 
unitaries $w_g\in {\Cal U}(N_2)$ for all $g\in G_1$ such that $(i), (ii)$  
are satisfied, then $\Phi$ can be extended to an isomorphism 
$N_1\rtimes_{(\sigma_1, v_1)} \Gamma_1\simeq N_2\rtimes_{(\sigma_2, v_2)} \Gamma_2$ $($hence, to 
an isomorphism of the associated inclusions$)$.

\endproclaim

\vskip.05in

\noindent
{\bf 1.3.   1-cocycles for actions.} Assume $\sigma$ is a 
genuine action of $\Gamma$ on the II$_1$ factor $N$. 
A map $w : \Gamma \rightarrow \Cal U(N)$ satisfying condition 
$$
w_g \sigma_g(w_h) = w_{gh}, \ \forall g, h \tag 1.3.1
$$

\noindent is called a 1-{\it cocycle for} $\sigma$.  Such a 1-cocycle for 
$\sigma$ is a {\it coboundary} (or it is {\it trivial}) if there exists 
a unitary element $v \in \Cal U(N)$ such that $w_g = v^* \sigma_g(v), 
\ \forall g$. (Clearly, such maps $w_g$ do satisfy the 1-cocycle 
condition (1.3.1)). 

The map $w$ is called a {\it weak} 1-{\it cocycle} if 
it satisfies the relation (1.2.1) modulo the scalars, 
i.e., 
$$
w_g \sigma_g(w_h) w_{gh}^* \in \Bbb T 1, 
\forall g,h \in \Gamma \tag 1.3.1'
$$ 
Note that this is equivalent to Ad$(w_g)\sigma_g$ being an action. 
Note also that if $w$ is a weak 1-cocycle then 
$\mu_{g,h} = 
w_g \sigma_g(w_h) 
w_{gh}^*$ is a scalar 2-cocycle for $\Gamma$, i.e., $\mu\in H^2(\Gamma)$. 
Also, cocycle conjugacy of two (genuine) actions $\sigma_i : \Gamma \rightarrow \text{\rm Aut}(N_i)$, $i=1,2$, 
amounts to conjugacy of $\sigma_1$ and $\sigma'_2$, where $\sigma'_2(g)=\text{\rm Ad}(w_g) \sigma_2(g)$, $g\in \Gamma$, 
for some $1$-cocycle $w$ for $\sigma_2$.

A (weak) 
1-cocycle $w$ is {\it weakly trivial}  
(or {\it weak cobouboundary}) if there exists 
a unitary element $v \in \Cal U(N)$ such that $v w_g \sigma_g(v)^*  \in 
\Bbb T 1,  \forall g$. 

Two (weak) 1-cocycles $w, w'$ of the action $\sigma$ are 
equivalent if there exists a unitary element 
$v\in N$ such that $w'_g = vw_g\sigma_g(v)^*, \forall g\in \Gamma$ 
(resp. modulo scalars). Thus, a weak 1-cocycle is weakly trivial 
iff it is equivalent to a scalar valued weak 1-cocycle (N.B.: these are 
just plain scalar functions on $\Gamma$). 
Note that the scalar valued genuine 1-cocycles 
are just characters of 
$\Gamma$. 

Two free actions $\sigma_1, \sigma_2$ of $\Gamma$ on $N$ 
are cocycle conjugate iff $\sigma_1$ is conjugate to $\sigma_2'$, where $\sigma_2'(g)=\text{\rm Ad}(w_g)\sigma_2(g)$,  $\forall g\in \Gamma$, 
for some $1$-cocycle $w$ for $\sigma_2$. 

We also mention here a well known result from [J80], showing that any 1-cocycle of an 
action of a finite group $\Gamma$ is co-boundary. This property is actually specific to finite groups: we use a result 
in [P01a] to deduce that if $\Gamma$ is infinite, then there exist free ergodic actions $\Gamma \curvearrowright R$ 
which admit non-trivial 1-cocycles. (N.B. In the particular case when $\Gamma$ is amenable, this fact can be derived from [Oc]  as well).  

\proclaim{Proposition} $1^\circ$ Let $\Gamma \curvearrowright^\sigma N$ be a free action of a finite group $\Gamma$ on a 
$\text{\rm II}_1$ factor $N$. If $w$ is a $1$-cocyle for $\sigma$, then there exists $u\in \Cal U(N)$ such that $w_g=u\sigma_g(u^*)$, $\forall g\in \Gamma$. 

$2^\circ$ Let  $(N_0, \varphi_0)$ be a copy of the $2$ by $2$ matrix algebra with the state given by weights $\{\frac{1}{1+\lambda}, \frac{\lambda}{1+\lambda}\}$,  
for some $0<\lambda<1$, and $\Gamma$ be an infinite group. Let $\Gamma \curvearrowright (\Cal N, \varphi)=\overline{\otimes}_g (N_0, \varphi_0)_g$ 
be the  Bernoulli $\Gamma$-action with base $(N_0, \varphi_0)$. Let $\Gamma \curvearrowright^\sigma N=\Cal N_\varphi\simeq R$ be 
the corresponding Connes-St\o rmer Bernoulli action. Let $B\subset N$ be an atomic von Neumann subalgebra of the form $\oplus_n B_n$,  
with $B_n\simeq M_{k_n \times k_n}(\Bbb C)$ having minimal projections   
of trace $\lambda^{m_n}$, with $m_1< m_2 < ....$.   
Then there exists a $1$-cocycle $w$ for $\sigma$ 
such that $\sigma'_g=\text{\rm Ad}(w_g)\sigma_g$, $g\in \Gamma$, has $B$ as its fixed point algebra. 
If $B\neq \Bbb C$, then any such $1$-cocycle is not a co-boundary. 
\endproclaim
\noindent
{\it Proof}. $1^\circ$ This is (Corollary 2.16 in [CT76; see also [J80]).  
We include here the proof, for completeness, which is based on Connes well known ``2 by 2 matrix trick''. 

Thus, let
$\tilde{\sigma}$ be the action of $\Gamma$ on $\tilde{N}=M_{2 \times
2}(N)=N \otimes M_{2 \times 2}(\Bbb C)$ given by
$\tilde{\sigma}_g=\sigma_g \otimes id$. If $\{e_{ij} \mid 1\leq i,j
\leq 2\}$ is a matrix unit for $M_{2 \times 2} (\Bbb C)\subset
\tilde{N}$, then $\tilde{w}_g=e_{11}+ w_g e_{22}$ is a cocycle
for $ \tilde{\sigma}$. If $Q\subset \tilde{N}$ denotes the fixed
point algebra of the action $\tilde{\sigma}_g'=\text{\rm Ad}(\tilde{w}_g) \tilde{\sigma}$,
then $e_{11}, e_{22}\in Q$ and the existence of a unitary element
$u\in N$ satisfying $w_g=u\sigma_g(u^*)$, $\forall g$, is equivalent to the fact that $e_{11},
e_{22}$ are equivalent projections in $Q$. But the fixed point algebra of any free action of a finite group on a II$_1$ factor is a II$_1$ factor. 
Thus, $e_{11}, e_{22}$ are equivalent in $Q$ and $w$ follows co-boundary. 

$2^\circ$ For each $n\geq 1$, let  $\{V^j_n\}_{1\leq j \leq k_n} \in \Cal N$ 
be isometries such that $V^j_n{V^j_n}^*\in N$, $\tau(V^j_n{V^j_n}^*)=\lambda^{m_n}$, $V^j_nN{V^j_n}^* = N$ and $\{V^i_n{V^j_n}^* \mid 1\leq i,j \leq k_n\}$ 
be the matrix units of $B_n$.  
Let also $\pi_n$ be the trivial representation of $\Gamma$ of multiplicity $k_n$. Then by 
(Theorem 3.2 in [P01a]), $w_g =\sum_n \sum_{i,j}\pi_n(g)_{i,j} V^i_n \sigma_g(V^j_n)^*= \sum_n \sum_i V^i_n \sigma_g(V^i_n)^*$, $g\in \Gamma$, 
defines a 1-cocycle for $\sigma$ and $\sigma'_g=\text{\rm Ad}(w_g) \sigma_g$ has $B$ as its fixed point algebra. 

Since the fixed point algebra of an action is a conjugacy invariant of the action and $\sigma$ is mixing (thus ergodic), it follows that $\sigma, \sigma'$ are 
not conjugate, in particular there exists no $u\in \Cal U(N)$ such that $ \sigma'_g=\text{\rm Ad}(u)\sigma_g\text{\rm Ad}(u^*)$, $\forall g$, a relation 
that amounts to $w_g=u\sigma_g(u^*)$ modulo scalars, $\forall g$.    
\hfill $\square$

\vskip.05in
\noindent
{\bf 1.4. Vanishing cohomology  and property  $\Cal V\Cal C$.} 
The 2-cocycle $v$ for the cocycle action $\sigma$ {\it vanishes} (or it is a 
{\it coboundary}) if there exists a map $w : \Gamma \rightarrow \Cal U(N)$ such 
that $w_e=1$ and $v = \partial w$, i.e.:
$$
v_{g, h} = (\partial w)_{g,h} \overset \text{\rm def} \to = \sigma_g(w^*_h) 
w^*_g w_{gh}, \ \forall g, h \in \Gamma. \tag 1.4.1
$$

The 2-cocycle $v$ {\it weakly-vanishes} (or it is a {\it weak coboundary}) if there exists $w : \Gamma  
\rightarrow \Cal U(N)$ such that $w_e=1$ 
and $v = \partial w$ modulo scalars, i.e.:
$$
w_g \sigma_g(w_h)v_{g,h} w_{gh}^* \in \Bbb C1, \ 
\forall g, h \in \Gamma. \tag 1.4.2
$$

Note that this is equivalent to
$$
\left( \text{\rm Ad} \ w_g \sigma_g \right) 
\left( \text{\rm Ad} \ w_h \sigma_h 
\right) = \ \text{\rm Ad} \ w_{gh} \sigma_{gh}, \ \forall g, h 
\tag 1.4.2$^\prime$
$$
i.e., to $\sigma'_g \overset \text{\rm def} \to = \ \text{\rm Ad} \ w_g 
\sigma_g$ being a ``genuine'' action. 

In turn, the vanishing of $v$ is equivalent to the existence of unitary elements $\{w_g\}_g \subset N$ 
such that $U'_g = w_gU_g\in M=N\rtimes_{(\sigma, v)} \Gamma$ give a representation of $\Gamma$ 
(i.e., $U'_gU'_h=U'_{gh}, \forall g, h\in \Gamma$). 

Given a II$_1$ factor $N$, we denote by $\Cal V\Cal C(N)$ the class of countable groups $\Gamma$ for which 
any free cocycle action $(\sigma, v)$ of $\Gamma$ on $N$ has the property that the 
$2$-cocycle $v$ vanishes (or is coboundary) and 
by $\Cal V\Cal C_w(N)$ the class of groups $\Gamma$ for which any free 
cocycle action $(\sigma, v)$ of $\Gamma$ on $N$ has the property that $v$ is a weak-coboundary. 

We denote by $\Cal V\Cal C$ (respectively $\Cal V\Cal C_w$) 
the class of countable groups $\Gamma$ with the property that $\Gamma \in \Cal V\Cal C(N)$ (resp. $\Gamma \in \Cal V\Cal C_w(N)$) 
for any II$_1$ factor $N$.  If $\Gamma\in \Cal V\Cal C$ then we also say that $\Gamma$ has {\it property} $\Cal V\Cal C$ or that it is a $\Cal V\Cal C$-{\it group}. 

We are especially interested in identifying the $\Cal V\Cal C$ and $\Cal V\Cal C_w$ groups, i.e.,   groups that 
have the most ``universal'' 
vanishing cohomology property. Other classes of interest  will be $\Cal V\Cal C(N)$ for $N$ equal to the hyperfinite II$_1$ factor $R$  
and for $N$ equal to the free group factor $L(\Bbb F_\infty)$. This is because $R$ and $L(\Bbb F_\infty)$ are the most interesting 
``non-commutative probability spaces''.  Also, any countable group $\Gamma$ has ``many'' 
free actions on these factors, in fact both of them have a lot of generalized symmetries (notably $L(\Bbb F_\infty)$, on which by [PS01] 
any ``group-like'' object admits  free actions). Also, both factors admit many cocycle actions that are ``hard to untwist'' (cf. [CJ84], [P01a] and Section 3 and 4 below). 
(N.B. It should be noticed that by the way we have defined  
$\Cal V\Cal C(N)$, if a factor $N$ has only inner automorphisms, i.e., Out$(N)=\{1\}$, like the examples in [IPeP05], then any $\Gamma$ belongs to 
$\Cal V\Cal C(N)$!)

We'll now show that the class $\Cal V\Cal C$ is closed to amalgamated free products over finite subgroups and that 
vanishing cohomology for cocycle actions of countable groups is  essentially a ``separability'' property: 

\proclaim{1.5. Proposition} $1^\circ$ if $\{\Gamma_n\}_{n\geq 0} \subset \Cal V\Cal C(N)$ $($respectively $\Cal V\Cal C_w(N))$ 
for some $\text{\rm II}_1$ factor $N$ and $K\subset \Gamma_n$ is a common finite subgroup, $n\geq 0$, 
then $\Gamma_0*_K \Gamma_1*_K \Gamma_2 *_K ...\in \Cal V\Cal C(N)$ $($respectively $\Cal V\Cal C_w(N))$. 
Also, if $\{\Gamma_n\}_n \subset \Cal V\Cal C$ 
$($resp. $\Cal V\Cal C_w)$, then $\Gamma_0*_K \Gamma_1 *_K ....  \in \Cal V\Cal C$  $($resp. $\Cal V\Cal C_w)$. 

$2^\circ$ Let $N$ be a $\text{\rm II}_1$ factor, $\Gamma \subset \text{\rm Out}(N)$ a countable group with a lifting 
$\{\sigma_g\}_{g\in \Gamma}\subset  \text{\rm Aut}(N)$ and denote $v_{g,h}\in \Cal U(N)$ a set of unitaries 
satisfying $\text{Ad}(v_{g,h})=\sigma_g\sigma_h\sigma_{gh}^{-1}$, $g, h\in \Gamma$.  
There exists a separable II$_1$ subfactor $Q\subset N$ that contains 
the countable set $\{v_{g,h} \mid g, h\in \Gamma\}$ and is normalized by $\sigma$, with $\sigma_g$ 
outer, $\forall g\in \Gamma$. 

$3^\circ$ $\Cal V\Cal C$ $($respectively $\Cal V\Cal C_w)$ coincides with 
the class of countable groups $\Gamma$ with the property that $\Gamma \in \Cal V\Cal C(N)$ $($resp. $\Gamma \in \Cal V\Cal C_w(N))$  
for any separable $\text{\rm II}_1$ factor $N$.  

\endproclaim
\noindent
{\it Proof}. $1^\circ$ Assume $\Gamma_n \in \Cal V\Cal C(N)$. 
Let $(\sigma, v)$ be a free cocycle action of $G=*_K \Gamma_n$ on $N$ and denote 
$M=N \rtimes_\sigma G$ with $U_g, g\in G$ the corresponding canonical unitaries. Since  
$\Gamma_n \in \Cal V\Cal C(N)$, there exist 
unitaries $\{w^n_g \mid g \in \Gamma_n\}$   in $N$ such that $U^n_g=w^n_g U_g, g\in \Gamma_n$, 
give left regular representations of $\Gamma_n$. Replacing $U_g$ by $w^0_gU^0_g$, $g\in \Gamma_0$, 
we may assume $w^0_g=1$, $\forall g\in \Gamma_0$. 

But then for each $n\geq 1$, $U^n_k=w^n_kU_k, k\in K$, for some 
1-cocyles $w^n: K \rightarrow \Cal U(N)$ for the restriction to $K$ of the 
$\Gamma_n$-action $\sigma_n$ implemented by $U^n_g, g\in \Gamma_n$. 
By ([J80]; see Proposition 1.3 above) any  1-cocycle of a free action of a finite group vanishes.  
Hence, there exists $v_n\in \Cal U(N)$ (with $v_0=1$) such that $w^n_k = v_n\sigma_n(k)(v_n^*)$, equivalently $U^n_k=v_nU_k v_n^*$, $k\in K$. 
But then the unitaries $\{v_n^*U^n_g\sigma_g(v_n)\mid g\in \Gamma_n, n\geq 0 \}$ generate inside $M$ a copy of the left 
regular representation of $G=*_K \Gamma_n$ implementing a $G$-action on $N$ that gives an inner perturbation 
of the initial cocycle $G$-action  $\sigma$. 

$2^\circ$ We construct recursively an increasing sequence of separable 
von Neumann subalgebras $Q_n, n\geq 0$, such that $Q_0\supset \{v_{g,h} \mid g, h \in \Gamma\}$ and for each 
$m\geq 1$ we have 

\vskip.05in

$(a)$ $E_{Q_m'\cap N}(x)=\tau(x)1, \forall x\in (Q_{m-1})_1$;  

$(b)$ $E_{Q_m'\cap N\rtimes \Gamma}(U_g)=0, \forall g\neq e$; 

$(c)$ $Q_m \supset \cup_g \sigma_g(Q_{m-1})$,

\vskip.05in 
\noindent
where $U_g\in N \rtimes \Gamma$ are the canonical unitaries implementing $\sigma$. 

Assume we have constructed these algebras up to $m=n$. 
Since $N'\cap N\rtimes \Gamma = \Bbb C$,  by using (Theorem 0.1 in [P13]) we can get a Haar unitary $v=(v_k)_k \in N^\omega$ that's free 
independent to $Q_{n-1} \cup \{U_g\}_g$. Thus, if we take $Q^0_n$ to be the von Neumann algebra generated by 
$Q_{n-1}$ and $\{v_k\}_k$, then we already have $(a)$ and $(b)$ satisfied for $Q_n=Q_n^0$, and then we can replace this ``initial'' $Q^0_n$ 
by the von Neumann algebra generated by $Q_n=\cup_g \sigma_g(Q^0_n)$, to have $(c)$ satisfied as well.  

Finally, if we define $Q=\overline{\cup_n Q_n}^w$, then $Q$ is clearly separable, condition $(c)$ insures that 
$Q$ is normalized by $\sigma$, condition $(a)$ shows that $E_{Q'\cap N}(\cup_n Q_n)\in \Bbb C1$, implying that $Q$ 
is a factor, while condition $(b)$ shows that $\sigma_g$ is outer on $Q$, $\forall g\neq e$. 

$3^\circ$ This part is now trivial by $2^\circ$. 

\hfill $\square$

\vskip.05in 

\noindent
{\bf 1.6. Remarks} $1^\circ$ As we will see in Sections 3 and 4, it is in general not true that 
if $\Gamma_i$ are in $\Cal V\Cal C$ then their amalgamated free product over a common (infinite) 
amenable subgroup $H\subset \Gamma_i$, $\Gamma=\Gamma_1 *_H \Gamma_2$ is in $\Cal V\Cal C$. 
For instance, $\Bbb Z^2 \rtimes SL(2, \Bbb Z)$ does not even belong to $\Cal V\Cal C_w(R)$ (see Theorem 3.2).  

$2^\circ$ The classes $\Cal V\Cal C$ may satisfy other general permanence properties.  For instance, 
it may be true that $\Gamma\in \Cal V\Cal C$ 
implies $\Gamma_0 \in \Cal V\Cal C$ for any  subgroup $\Gamma_0\subset \Gamma$ 
(or at least when $[\Gamma:\Gamma_0]< \infty$). However, the obvious idea for a proof, which is to ``co-induce'' 
a given cocycle action $\Gamma_0 \curvearrowright^{\sigma_0} N_0$ to a set of automorphisms 
$\{\sigma_g \mid g\in \Gamma\}$ on $N=N_0^{\overline{\otimes}\Gamma/\Gamma_0}$ doesn't work when $[\Gamma:\Gamma_0]=\infty$, because 
an infinite tensor product of inner automorphisms may become outer, so the  
$\sigma_g$'s may in fact not give a cocycle action of $\Gamma$. When the index of $\Gamma_0$ 
in $\Gamma$ is finite, then $\sigma$ defined this way does give a cocycle action of $\Gamma$ on $N$, 
but it is not immediate of how to ``bring down to $N_0$'' the vanishing of the cohomology for $\sigma$ to 
the vanishing of the cohomology for the initial $\Gamma_0 \curvearrowright^{\sigma_0} N_0$.

\heading 2. Groups with the property $\Cal V\Cal C$   \endheading 

In this section we prove a vanishing cohomology result for arbitrary  free cocycle actions of countable 
amenable groups on arbitrary II$_1$ factors. 

To do this, we'll first show that any 
amenable subgroup $\Gamma \subset \text{\rm Aut}(N)/\text{\rm Inn}(N)$ can be lifted 
to a set $\{\sigma_g\mid g\in \Gamma\}\subset \text{\rm Aut}(N)$ normalizing a ``large'' hyperfinite subfactor of $N$ (see Theorem 2.1).
As it happens, this property, which is interesting by itself, characterizes 
the amenability of the group $\Gamma$. Indeed, we will show in Section 5 that any non-amenable group admits a free action 
on $N=L(\Bbb F_\infty)$ that cannot be perturbed to a cocycle action that normalizes a hyperfinite subfactor of $N$.  

Once we prove that 
any cocycle action $\sigma$ of an amenable group $\Gamma$ on $N$ 
normalizes (modulo inner perturbation) a hyperfinite II$_1$ factor $R\subset N$, 
we reduce the vanishing cohomology problem to the case $N=R$,   
where by a well known result of Ocneanu [Oc85] free cocycle actions  
of amenable groups 
can indeed be ``untwisted'' to genuine actions. The fact that $R$ is ``large in $N$'' assures   
that by untwisting $\sigma$ on $R$ we have  untwisted it as an action on $N$ as well.  

To show that $\sigma$ normalizes up to Inn$(N)$ a  ``large hyperfinite subfactor of $N$'', we reduce the problem to a statement about 
commuting squares of subfactors, as follows. Assume $\Gamma$ is generated 
by a finite set $e\in F=F^{-1}\subset \Gamma$ and consider the {\it locally trivial subfactor} obtained by the diagonal embedding 
$N^{\sigma,F}:=\{\sum_{g\in F} \sigma_{g}(x) e_{gg} \mid x\in N\} \subset M_{|F| \times |F|}(N)=:M^{\sigma,F}$, where $\{e_{gh}\}_{g,h\in F} 
\subset M_{|F| \times |F|}(\Bbb C)$ are the matrix units (see 5.1.5 in [P91]). If $Q\subset \Cal R$ is an inclusion of factors with 
$Q\subset N^{\sigma,F}$, $\Cal R \subset M^{\sigma,F}$, 
$e_{gg}\in \Cal R, \forall g\in F$,  and $(Q\subset \Cal R) \subset (N^{\sigma,F} \subset M^{\sigma,F})$ makes a non-degenerate commuting square,  
then $Q\subset \Cal R$ is itself locally trivial and there exist unitary elements $w_g \in N$ such that $w_ge_{eg}\in \Cal R$.  If one denotes $Q_0\subset N$ the image  
of $Q$ under the isomorphism $N^{\sigma,F} \simeq N$, then this amounts to $Q_0$ being invariant to $\sigma'_{g}=\text{\rm Ad}w_g\circ \sigma_{g}$, 
$\forall g\in F$. Moreover, if $Q'\cap \Cal R={N^{\sigma,F}}' \cap \Cal R$, then ${\sigma'_g}_{|Q_0}$ is outer iff $\sigma_{g}$ 
is outer, $\forall g\in F$.  This observation applied to $N^{\sigma,F} \subset M^{\sigma,F}_n$  
(where $M^{\sigma,F}_n$ are the factors in the tower for $N^{\sigma,F} \subset M^{\sigma,F}$), 
in combination with (5.1.5 in [P91]), shows that if all the higher relative commutants in the Jones towers for 
$Q\subset \Cal R$ and $N^{\sigma, F} \subset M^{\sigma,F}$ coincide, then $\sigma'$ implements an outer action of $\Gamma$ on $Q_0$. 

So all we need to do is to produce an 
inclusion of hyperfinite factors $Q\subset \Cal R$ inside $N^{\sigma,F} \subset M^{\sigma,F}$, making a commuting square and having same higher relative 
commutants.  

We will solve  this commuting square problem by only using that $N^{\sigma,F} \subset M^{\sigma,F}$ has amenable graph. 
Thus, we will in fact prove that any finite index 
inclusion of separable II$_1$ factors $N\subset M$  with amenable standard invariant $\Cal G_{N\subset M}$ 
contains an inclusion of hyperfinite factors $(Q\subset R)\subset (N\subset M)$, 
that makes a non-degenerate commuting square and has identical 
higher relative commutants in the associated Jones tower (in particular same standard invariant),  
in fact even satisfies  the  {\it strong smoothness} condition 
$Q'\cap R_n=Q'\cap M_n= N' \cap M_n$, $\forall n$ (see Theorem 2.10 below).  We'll obtain $Q\subset R$ as an inductive limit of relative commutants 
$P_n'\cap N \subset P_n'\cap M$ of  a decreasing sequence of finite index subfactors $P_n\subset N$ that come from repeated downward basic constructions of  subfactors 
$Mp'\subset p'M_np'$ obtained 
by appropriate induction/reduction in the Jones tower $N\subset M\subset M_1 \subset ...$, with choices ``dictated'' by the 
local characterization of the amenability of $\Gamma_{N\subset M}$ in ([P97], Theorem 6.1).

\proclaim{2.1. Theorem} Let $N$ be a  $\text{\rm II}_1$ factor and $\sigma:\Gamma \rightarrow \text{\rm Aut}(N)$ an outer 
action of a countable amenable group $\Gamma$  on $N$, with $H^3(\Gamma)$-obstruction $\nu_\sigma$ and with $v_{g,h}\in \Cal U(N)$ 
satisfying $\sigma_g\sigma_h = \text{\rm Ad}(v_{g,h})\sigma_{gh}$, $\forall g, h\in \Gamma$. 
Then there exist $\{w_g\}_g \subset \Cal U(N)$ and a hyperfinite subfactor $R\subset N$ 
such that if we denote $\sigma'_g=\text{\rm Ad}(w_g) \sigma_g$ and $v'_{g,h}= w_g\sigma_g(w_h)v_{g,h}w_{gh}^*$, $g, h\in \Gamma$, 
then we have: 

\vskip.05in
$(2.1.1)$ $\sigma'_{g}(R)=R$ and $v'_{g,h}\in R$, $\forall g, h\in \Gamma$. 

\vskip.05in

$(2.1.2)$ $\{{\sigma'_g}_{|R}\}_g$ is an outer action of $\Gamma$ on $R$ with same $H^3(\Gamma)$-obstruction as $\sigma$. 

\vskip.05in
If in addition $(\sigma, v)$ is a free cocycle action of $\Gamma$ on $N$, then $\{w_g\}_g$ can 
be chosen so that $(\sigma'_{|R}, v')$ is a free cocycle action of $\Gamma$ on $R$. 
Moreover, if $N$ is separable, then one can choose $\sigma',  v', R$ so that to also have $R'\cap N=\Bbb C$. 
\endproclaim

Let us show right away how Theorem 2.1 combined with Ocneanu's Theorem in [Oc85] can be used to derive the 
vanishing cohomology result for cocycle actions of arbitrary amenable groups: 

\newpage

\proclaim{2.2. Theorem} Let $N$ be a $\text{\rm II}_1$ factor, $\Gamma$ a countable amenable group 
and $(\sigma, v)$  a free cocycle action of $\Gamma$ on $N$. Then there exist  unitary elements 
$\{w_g\in \Cal U(N) \mid g\in \Gamma\}$ such that $\sigma'_g=\text{\rm Ad}(w_g)\circ \sigma_g$, $g\in \Gamma$, is a 
genuine action of $\Gamma$ on $N$, in fact such that moreover we have $v_{g,h}=\sigma_g(w_h^*)w_g^*w_{gh}$, $\forall g, h \in \Gamma$. 
\endproclaim
\noindent
{\it Proof  of} 2.2. By Theorem 2.1, there exist unitary elements $w^0_g\in N$, $g\in \Gamma$, 
and a hyperfinite II$_1$ subfactor $R\subset N$ such that: 

\vskip.05in 
$(2.2.1)$ $\sigma^0_g:=\text{\rm Ad}(w^0_g) \sigma_g$ leaves $R$ invariant, $\forall g\in \Gamma$; 

\vskip.05in

$(2.2.2)$  $v^0_{g, h}:= w_g\sigma_g(w_h)v_{g,h}w_{gh}^*$ belongs to $R$, $\forall g, h\in \Gamma$.   

\vskip.05in

$(2.2.3)$ $\sigma^0_{g|R}$ is outer, $\forall g\neq e$. 

\vskip.05in

But then, $(\sigma^0_{|R}, v^0)$ implements a free cocycle action of the countable amenable group $\Gamma$ on $R$, 
so by Ocneanu's Theorem [Oc85] the 2-cocycle $v^0$ is co-boundary on $R$,  i.e., 
there exist unitary elements  $w^1_g\in R$ such that 
$v^0_{g,h}=\sigma^0_g({w^1_h}^*){w^1_g}^*w^1_{gh}$, $\forall g, h \in \Gamma$. This shows that $w_g=w^1_gw^0_g$ 
satisfy the required condition. 
\hfill $\square$

\proclaim{2.3. Corollary} The class $\Cal V\Cal C$ contains all countable amenable groups and is closed to free products with 
amalgamation over finite subgroups, i.e., if $\{\Gamma_n\}_n \subset \Cal V\Cal C$ and $K\subset \Gamma_n$ is a common finite subgroup, 
then $\Gamma_0*_K \Gamma_1 *_K ....  \in \Cal V\Cal C$. 

\endproclaim

For the rest of this section, we will use concepts and notations from  [J83] (such as 
the basic construction, the Jones tower of factors, etc), as well as from ([PiP84], [P91], [P93], [94a], [94b], [P97]).  In particular, we will often use 
as framework 
the {\it symmetric enveloping} (abbreviated {\it SE}) inclusion $M\overline{\otimes} M^{^{op}}\subset M\underset{e_N}\to{\boxtimes} M^{^{op}}$ of $N\subset M$,   
introduced in [P94b] (cf. also the extended version of the paper, [P97]).  

We begin by recalling some properties relating Jones tower/tunnel of a subfactor with its symmetric enveloping inclusion. It will be useful for the reader to keep 
in mind that if $M \subset M_1 \subset ...$ is the Jones tower of factors associated with a 
subfactor of finite index $N\subset M$,  then $_ML^2(M_n)_M 
=_ML^2(M_1)\overline{\otimes}_M .... \overline{\otimes}_ML^2(M_1)_M$ ($n$-times tensor$/M$ product). Also, if one denotes 
by $\{\Cal H_k\}_{k\in K}$ the list of irreducible Hilbert $M$-bimodules appearing in $\{L^2(M_n)\}_n$, then for any $m\geq 1$ 
and any nonzero projection $p'\in M_1'\cap M_m$ 
we have $_ML^2(M_1)_M \subset _ML^2(p'M_mp')_M$, and thus we recover all $\{\Cal H_k\}_k$ in the tensor powers $_ML^2(p'M_mp')^{n\overline{\otimes}_M}_M$, $n\geq 1$. 
Equivalently, if $M \hookrightarrow p'M_mp' =M^0$ and $M^0\subset M^0_1 \subset ...$ is its Jones tower, then 
$\{\Cal H_k\}_k$  coincides with the list of irreducible Hilbert $M$-bimodules appearing in $L^2(M^0_n), n\geq 1$ (see 
[P91] and [Bi97] for basics of subfactor theory).

\proclaim{2.4. Lemma} Let $N\subset M$ be an extremal inclusion of $\text{\rm II}_1$ factors 
of index $[M:N]=\lambda^{-1}< \infty$,  
$T=M\overline{\otimes} M^{^{op}}\subset M\underset{e_N}\to{\boxtimes} M^{^{op}}=S$ 
its $\text{\rm SE}$ inclusion of $\text{\rm II}_1$ factors and 
$... \subset N_m \subset ... \subset N\subset^{e_{-1}} M \subset^{e_N=e_0} M_1 \subset^{e_1} ... M_m\subset^{e_m} .... $  a tunnel-tower for 
$N\subset M$ inside $S$. 

$1^\circ$ If $e^n_{-n}\in M_{n+1}\subset S$ denotes the projection of trace $[M:N]^{-n-1}$ obtained as a scalar multiple of the word 
of maximal length in $e_{-n}, ..., e_0, ...., e_n$, then $e^n_{-n}$ implements $E^M_{N_n}, E^{M^{^{op}}}_{N_n^{^{op}}}$ 
and one has $vN(T, e^n_{-n})=S$.  Thus, the map that acts as the identity on $M\vee M^{^{op}}$ and takes $e_{N_n}$ to $e^n_{-n}$ 
gives a natural identification between $M\vee M^{^{op}}\subset M\underset{e_{N_n}}\to{\boxtimes} M^{^{op}}$ and $T \subset S$. 

$2^\circ$ Let $p\in \Cal P(N_n'\cap M)$ and $p^{^{op}}\in M'\cap M_{n+1}\subset S$ 
its image under the antiautomorphism $^{^{op}}$ of $S$. Then $N_{n+1}pp^{^{op}}\subset pMpp^{^{op}} \subset pp^{^{op}}M_{n+1}pp^{^{op}}$ 
is a basic construction 
for $(V\subset U)=(N_{n+1}pp^{^{op}}\subset pMpp^{^{op}})$, with Jones projection $f=\tau(p)^{-1} pp^{^{op}}e^n_{-n}pp^{^{op}}=
\tau(p)^{-1}pe^n_{-n}p=\tau(p)^{-1} p^{^{op}}fp^{^{op}}$. Moreover $U\underset{e_V}\to{\boxtimes}U^{^{op}}$ is naturally 
embedded into $S$ as the  von Neumann subalgebra 
generated by $pMpp^{^{op}}$, $pp^{^{op}}M^{^{op}}p^{^{op}}$ and $f$. 
If in addition $p\in N_n'\cap N$, then this latter algebra is actually equal to $pp^{^{op}}Spp^{^{op}}$, thus giving a natural identification between  
$(U \overline{\otimes} U^{^{op}}\subset U\underset{e_V}\to{\boxtimes}U^{^{op}})$ and the amplification 
by $\tau(p)^2$ of $(T\subset S)$, with $e_V\mapsto f$.

$3^\circ$ Let $p\in \Cal P(N_n'\cap N)$ be as in the last part of $2^\circ$. Let $P\subset N$ be a subfactor such that $P\subset M$ 
is a downward basic construction for $M\simeq Mp^{^{op}} \subset p^{^{op}}M_{n+1}p^{^{op}}$ and denote 
$M\overline{\otimes} M^{^{op}}=T\subset S_1=M\underset{e_P}\to{\boxtimes} M^{^{op}}$ its $\text{\rm SE}$ inclusion. If $\{m_j\}_j\subset N$ 
is an orthonormal basis of $N$ over $P$, then $e=\sum_j m_je_Pm_j^*=\sum_j  {m_j^{^{op}}}^* e_P m^{^{op}}_j $
is a projection of trace $\lambda=[M:N]^{-1}$ in $S_1$ that 
implements  both $E^M_N$ and $E^{M^{^{op}}}_{N^{^{op}}}$ and satisfies $vN(T, e)=S_1$, thus giving an identification  
between $T\subset S$ and $T\subset S_1$, with $e_N \mapsto e$. 
\endproclaim
\noindent
{\it Proof}. Part $1^\circ$ is (Proposition $2.9 (a)$ in [P97]), the first part of $2^\circ$ is  
(Lemma $2.8. (c)$ and $2.9 (c)$ in [P97]), while the first part of $3^\circ$ is (Proposition 2.10 in [P97]). 

To prove the last part of $2^\circ$, note that with the notation $U=pMpp^{^{op}}$ and $S_0=U\underset{e_V}\to{\boxtimes} U^{^{op}}\subset pp^{^{op}}Spp^{^{op}}$ 
we have $_UL^2(S_0)_U\subset {_U}L^2(pp^{^{op}}Spp^{^{op}})_U$ as Hilbert bimodules. Then notice that   
by (Theorem 4.5 in [P97]), $_UL^2(pp^{^{op}}Spp^{^{op}})_U= \oplus_{k\in K} \Cal H'_k \overline{\otimes} \overline{\Cal H'}_k^{^{op}}$, 
where $\{\Cal H'_k\}_{k\in K}$ denotes the reduction by $pp^{^{op}}$ 
of  all distinct irreducible Hilbert $M$-bimodules in $\cup_n L^2(M_n)$. Since 
the list of irreducible $U$-bimodules in the Jones tower for $N_np\subset pMp$ contains all $\{\Cal H'_k\}_{k\in K}$ (because 
$p\in N_n'\cap N$), 
it follows that $_UL^2(pp^{^{op}}Spp^{^{op}})_U \subset   {_U}L^2(S_0)_U$ as well. Thus, the inclusion $_UL^2(S_0)_U$ 
$\subset {_U}L^2(pp^{^{op}}Spp^{^{op}})_U$ 
is an equality, forcing $S_0=S$ as well. 

The last part of $3^\circ$ follows by taking again into account (Lemma 2.8. (c) in [P97]), part $2^\circ$ above and the remark before the statement of the lemma,  
then using the same ``exhaustion by bimodules'' argument used above (based on Theorem 4.5 in [P97]). 
\hfill $\square$

\vskip.05in
\noindent
{\bf 2.5. Definition}. Let $N\subset M$ be an inclusion of II$_1$ factors with finite index and denote by  
$N\subset M \subset M_{1}\subset  M_2\subset...$ its Jones tower.  
A subfactor $P$ in $N$ is said to be $(N\subset M)$-{\it compatible}  
if there exist $n\geq 1$ and a non-zero projection $p'\in M_1'\cap M_n$ such that $Pp' \subset Mp'  
\subset p'M_np'$ is a basic construction. Let us note right away that this property is in some sense ``hereditary'':

\proclaim{2.6. Lemma} With $N\subset M$ as above, assume $P\subset N$ is $(N\subset M)$-compatible. If a subfactor $Q\subset P$ is $(P\subset M)$-compatible, 
then $Q\subset N$ is $(N\subset M)$-compatible. 

\endproclaim
\noindent
{\it Proof}. Let $P \subset N \subset M \simeq Mp' \subset M_1p'\subset p'M_np'$ be a basic construction, for some $n\geq 1$ and a non-zero projection $p'\in M_1'\cap M_n$.  
Note that if we denote $V=Mp'\subset p'M_np'=U$ then, given any $m\geq 1$, its associated Jones tower of factors up to $m$, 
$V_1 \simeq P \subset M \simeq V\subset U \subset U_1 \subset ... U_m$, can be 
realized (up to isomorphism) by inducing/reducing in the initial tower $M \subset M_1\subset M_2... \subset M_k$, for some large enough $k$, with the 
projections involved $p'_{i_j}\in M_{i_j}'\cap M_{i_{j+1}}$, where $i_0=1$, $p'_{i_0}=p'$, $i_1= n$, and $i_0 < i_1 < ...$.  
This shows that if $Q\subset P$ is $(P\subset M)$-compatible, then one obtains a basic construction $Q\subset M \simeq Mq'\subset q'M_{k_0}q'$, 
for some appropriate $k_0$ and $q'\in M_1'\cap M_{k_0}$ obtained as a product of such $p'_{i_j}$. But this means that  $Q\subset N$ 
is $(N\subset M)$-compatible. 
\hfill $\square$

\vskip.05in

For the reader's convenience, we recall here two of the equivalent definitions of amenability for ``group-like'' objects arising from subfactors, that we have  
introduced and studied in ([P91], [P93], [P94b], [P97]),  and that we need  hereafter.  

The {\it standard invariant} $\Cal G_{N\subset M}$ 
of an extremal inclusion of factors with finite Jones index $N\subset M$   is {\it amenable} if its  {\it principal graph}  $\Gamma_{N\subset M}$ 
satisfies the Kesten-type condition $\|\Gamma_{N\subset M}\|^2=[M:N]$.  

This very first definition of amenability was proposed in [P91],  and we will also take it to be the definition of amenability for 
the various abstractizations of standard invariants: a {\it standard $\lambda$-lattice} $\Cal G$ as in [P04b] (or a {\it planar algebra} as in [J99]) is {\it amenable}  
if its graph $\Gamma_{\Cal G}$ satisfies the condition $\|\Gamma_{\Cal G}\|^2=\lambda^{-1}$.

An alternative notion of amenability was introduced in [P93], by requiring the following F\o lner-type condition on $\Cal G=\Cal G_{N\subset M}$: 
let $\{d_k\}_{k\in K}$ denote the {\it standard weights} of $\Gamma_{N\subset M}$ (resp. $\Gamma_\Cal G$),  obtained for instance as 
$\text{\rm dim}(_M{\Cal H_k}_M)^{1/2}$, where 
$\{\Cal H_k\}_{k\in K}$ is the list of irreducible $M-M$ Hilbert bimodules appearing at even levels  in $\Cal G$, 
indexed by the set $K$ of left vertices of the bipartite 
graph $\Gamma_{\Cal G}=\Gamma_{N\subset M}$; 
$\Cal G$ satisfies the {\it F\o lner condition} if for any $\varepsilon >0$ there exists a finite set $F\subset K$ such that 
if one denotes by $\partial F=\{k\in K\setminus F \mid \exists k_0\in F$ with $(\Gamma_{\Cal G}\Gamma_{\Cal G}^t)_{kk_0}\neq 0\}$  
(the boundary of $F$) then $\underset{k\in \partial F}\to{\sum} d_k^2 < \varepsilon \underset{k\in F}\to{\sum} d_k^2$. 

These two conditions were shown equivalent in ([P97] Theorem 5.3; the result had already been announced in   
[P93] and  [P94b]).  Several other equivalent amenability conditions were in fact introduced and studied in [P97],  
notably a local finite dimensional approximation property ([P97], Theorem 6.1) which will be crucial for the proof of Theorem 2.10 below.  

The Kesten and the F\o lner-type conditions provide in particular equivalent definitions of amenability for a 
finitely generated rigid C$^*$-tensor category of Hilbert bimodules $\Cal G$ (as defined for instance in [PV14]), having 
$\{\Cal H_k\}_{k\in K}$ as irreducible morphisms. 
One then defines amenability for an arbitrary (possibly infinitely generated) rigid C$^*$-tensor category 
by requiring that any finitely generated subcategory is amenable (see the detailed definitions in the paragraphs preceding Theorem 2.12).

Due to its various equivalent characterizations, amenability in this context is a very ``robust'' property. For instance, since the SE inclusion 
$T=M\overline{\otimes} M^{^{op}}\subset M\underset{e_N}\to{\boxtimes} M^{^{op}}=S$ associated with an extremal inclusion $N\subset M$ 
coincides with the SE inclusion $M\overline{\otimes} M^{^{op}}\subset M\underset{e_{N_n}}\to{\boxtimes} M^{^{op}}$ associated 
with $N_n \subset M$, for any $n\geq 0$ (e.g., by Lemma 2.4.1$^\circ$ above), it follows from (Theorem 5.3 in [P97]) 
that $\Cal G_{N\subset M}$ is amenable iff $\Cal G_{M_i\subset M_j}$ is amenable for some (and thus all) $i< j$.  Note that 
this can also be deduced from the fact that $\Gamma_{M_i \subset M_j}$ is a alternate product of $\Gamma_{N\subset M}$ and its transpose, 
or of $\Gamma_{M\subset M_1}$ and its transpose, $j-i$ times, which shows that 
one always have $\|\Gamma_{M_i \subset M_j}\|=\|\Gamma_{N\subset M}\|^{j-i}$ (cf. 1.3.5 in [P91]), implying that $\|\Gamma_{N\subset M}\|^2=[M:N]$ 
iff $\|\Gamma_{M_i\subset M_j}\|^2=[M_j:M_i]$.  

Moreover, if $p'$ is a non-zero projection in $M_i'\cap M_j$, for some $i< j$ in $\Bbb Z$, then by (Coroally 6.6 $(ii)$ in [P97]) 
$\Gamma_{M_i\subset M_j}$ amenable (which 
we already know is equivalent to $\Gamma_{N\subset M}$ being amenable) 
implies $V=M_ip'\subset p'M_jp'=U$ has amenable graph as well. Also, note that if $j\geq i+1$ and $p'\in M_{i+1}'\cap M_j$, 
then by Lemma 2.4 and the above argument, one conversely has that $V\subset U$ amenable implies $N\subset M$ amenable. 

We record all these facts in the following:

\proclaim{2.7. Proposition}  Let $N\subset M$ be an extremal inclusion of factors with finite index and $\{M_i\}_{i\in \Bbb Z}$ be a tunnel/tower of factors for $N\subset M$. 
If $\Cal G_{N\subset M}$ is amenable, then for any $i<j$ in $\Bbb Z$ and any non-zero  projection $p'\in M_i'\cap M_j$, the inclusion $M_ip' \subset p'M_jp'$ has amenable standard invariant. If in addition $j\geq i+1$ and $p'\in M_{i+1}'\cap M_j$, then conversely $\Cal G_{M_ip'\subset p'M_jp'}$ amenable implies 
$\Cal G_{N\subset M}$ amenable. 
\endproclaim

\vskip.1in

\proclaim{2.8. Lemma} Let $N\subset M$ be a finite index extremal inclusion of $\text{\rm II}_1$ factors with  
amenable standard invariant and $\text{\rm SE}$ factor $S=M\underset{e_N}\to{\boxtimes} M^{^{op}}$. 
Given any $(N\subset M)$-compatible subfactor 
$P\subset N$ and any $\varepsilon>0$, there exists a $(P\subset M)$-compatible subfactor $Q\subset P$ 
such that $\|E_{(Q'\cap M)'\cap S}(x)-E_{M'\cap S}(x)\|_2 \leq \varepsilon \|x\|$, $\forall x\in P'\cap S$. 
\endproclaim
\noindent
{\it Proof}. Let us first note a few {\it Facts} needed in the proof. 

\vskip.05in
{\it Fact} 1.  It is sufficient to show that there exists a compatible subfactor 
$Q\subset P$ with the property that $\|E_{Q'\cap M)'\cap S}(f)-\tau(f)1 \|_2 \leq  \varepsilon/([M:P]+1)^{5/2}$, where 
$f\in S$ is the Jones projection for $P\subset M \subset \langle M, P\rangle$, viewed inside $S$ 
(cf. 2.4.$3^\circ$ above). 
\vskip.05in

Indeed, because if $\{m_j\}_j \subset M$ is an orthonormal basis of $M$ over $P$ 
with $[M:P]+1$ many elements of norm $\leq [M:P]^{1/2}$ (cf. [PiP84]), 
then $\{[M:P]^{1/2}m_j^{^{op}}f\}_j$ is an orthonormal basis of $P'\cap S$ over $M'\cap S=M^{^{op}}$ 
and any $x\in P'\cap S = \langle M^{^{op}}, f \rangle$ is of the form  $x=\sum_j [M:P]^{1/2} m_j^{^{op}}f y_j^{^{op}}$, 
where $y_j^{^{op}}=[M:P]^{1/2}E_{M^{^{op}}}(f{m_j^{^{op}}}^*x)\in M^{^{op}}$ has operator norm 
majorized by $[M:P]^{1/2} \|x\| \|m_jf\| = [M:P]^{1/2} \|x\|$, thus giving the estimates
$$
\|E_{(Q'\cap M)'\cap S}(x)-E_{M'\cap S}(x)\|_2 
$$
$$
=\|\sum_j [M:P]^{1/2} m_j^{^{op}} (E_{(Q'\cap M)'\cap S}(f)-\tau(f)1)y_j^{^{op}}\|_2 
$$
$$
\leq [M:P]^{1/2} \sum_j \|m_j\| \|y_j^{^{op}}\| \|E_{(Q'\cap M)'\cap S}(f)-\tau(f)1\|_2
$$
$$
\leq [M:P]^{3/2}([M:P]+1)\|x\|  \varepsilon/([M:P]+1)^{5/2} < \varepsilon \|x\|.
$$

\vskip.05in
{\it Fact} 2. By Proposition 2.7 above, $\Cal G_{N\subset M}$ amenable implies $\Cal G_{P\subset M}$ amenable. 

\vskip.05in

{\it Fact} 3. By (Theorem 6.1 in [P97]),  if $P\subset M$ is an extremal inclusion of factors 
with amenable standard invariant then for any $\delta >0$ there exists $n\geq 1$ and a projection $p\in P_n'\cap P$ 
such that $\|E_{{p(P_n'\cap M)p}'\cap p\langle M, P \rangle p}(f_0p)-\tau(f_0) p\|^2_2 
< \delta \tau(p)$, where 
$... \subset P_n\subset  P_{n-1}\subset ...P\subset M\subset^{f_0=e_P}  \langle M, P \rangle$ denotes a Jones tunnel 
and basic construction for $P\subset M$. 

\vskip.05in

Let us now proceed with the proof of the lemma. 
By Fact 2,  $P\subset M$ is amenable  so we can apply Fact 3 to $(P\subset M \subset^{f_0} \langle M, P \rangle)$,  
to get an $n\geq 1$ and a projection  $p\in P_n'\cap P$ such that 
$$
\|E_{{p(P_n'\cap M)p}'\cap p\langle M, P \rangle p}(f_0p)-\tau(f_0) p\|_2 
< \varepsilon \|p\|_2/([M:P]+1)^2  \tag 2.8.1
$$

By amplifying by $\alpha=\tau(p)^{-1}$ the inclusions of factors 

$$
P_np \subset pPp \subset pMp \subset^{f_0p} p\langle M, P\rangle p \tag 2.8.2
$$
using partial isometries in $P$, we obtain inclusions of factors 
$$
(P_np)^\alpha=Q\subset P \subset M \subset^{f_0} \langle M, P \rangle \tag 2.8.3
$$ 
having same relative commutants as $(2.8.2)$.  
Thus,  by $(2.8.1)$, it follows that $Q$ satisfies 
$$
\|E_{Q'\cap \langle M, P \rangle}(f_0)-\tau(f_0) \|_2 
< \varepsilon /([M:P]+1)^2,  \tag 2.8.4
$$

By the way it is defined, $Q\subset P$ is an $(P\subset M)$-compatible  
subfactor, and thus $(N\subset M)$-compatible as well, while by Fact 1  and $(2.8.4)$ we have 
$\|E_{(Q'\cap M)'\cap S}(x)-E_{M'\cap S}(x)\|_2 \leq \varepsilon \|x\|$, for all $x\in P'\cap S$.

\hfill $\square$

\proclaim{2.9. Lemma} Let $N\subset M$ be a finite index extremal inclusion of $\text{\rm II}_1$ factors with 
$N\subset M \subset M_1 .... \nearrow M_\infty$ its Jones tower of factors and 
$S$ its $\text{\rm SE}$ factor. If $B\subset M$ is a diffuse von Neumann subalgebra, then $B\not\prec_{M_\infty} M'\cap M_\infty$ 
and $B\not\prec_S M^{^{op}}$.  

\endproclaim
\noindent
{\it Proof}. By [P03], in order to prove $B\not\prec_{M_\infty} M'\cap M_\infty$, it is sufficient to prove that 
for any finite set $F$ in a given total subset  $X$ of $M_\infty$, there exist $u_n \in \Cal U(B)$ such that $\lim_n \|E_{M'\cap M_\infty}(y^*u_nx)\|_2=0$, 
$\forall x, y\in F$. Taking $X=\cup_m M_m$, it is sufficient to show this for any $m$ and any finite $F\subset M_m$. 
But by [P03] this amounts to $M\not\prec M'\cap M_m$, which is trivial since $B$ is diffuse and $M'\cap M_m$ is finite 
dimensional. 

To prove $B\not\prec_S M^{^{op}}$ we use the same criterion, but with $X=\cup_m M_m M^{^{op}}$, which is total in $S$ by (Section 4 in [P97]).   
Thus, if $F\subset X$ is finite then we may assume $F\subset M_mM^{^{op}}$ for some large $m$ so if $x=x_1x_2^{^{op}}, y=y_1y_2^{^{op}} \in F$,  
with $x_1, y_1 \in M_m, x_2, y_2\in M$, and we take $u_n \in \Cal U(B)$, then we get the estimate 

$$
\|E_{M^{^{op}}}({y_2^{^{op}}}^*{y_1}^*u_n x_1 x_2^{^{op}})\|_2
\leq \|x_2\| \|y_2\| \|E_{M'\cap M_m}(y_1^*u_n x_1)\|_2. 
$$

This shows that it is actually sufficient to check the criterion for $F\subset M_m$, which amounts again to $M\not\prec M'\cap M_m$ 
as before. 
\hfill $\square$

\proclaim{2.10. Theorem} Let  $N\subset M$ be a finite index extremal inclusion of separable $\text{\rm II}_1$ factors with  
amenable standard invariant. There exists a sub-inclusion of hyperfinite factors $(Q\subset R) \subset (N\subset M)$  
making a non-degenerate commuting square that's strongly smooth, i.e., $Q'\cap R_n=Q'\cap M_n= N' \cap M_n$ 
and $R'\cap R_n=R'\cap M_n=M'\cap M_n$, $\forall n$, 
where $N\subset M = M_0 \subset^{e_0} M_1\subset ...$  is the Jones tower for $N\subset M$ and $R_n=vN(R, e_0, ..., e_{n-1})$, $n\geq 1$, 
the tower for $Q\subset R$. Moreover, $Q\subset R$ can be obtained as 
$Q=\overline{\cup_n P_n'\cap N}\subset \overline{\cup_n P_n'\cap M}=R$, for some decreasing sequence of $(N\subset M)$-compatible 
subfactors $M\supset N \supset P_0\supset P_1 ...$. 
\endproclaim
\noindent
{\it Proof}. We split this proof into several parts. 

\vskip.05in

{\it Fact} 1. There exists a sequence of $(N\subset M)$-compatible subfactors $.... \subset P_n \subset ... P_0\subset N \subset M$ 
such that if we define $Q= \overline{\cup_m P_m'\cap N}\subset \overline{\cup_m P_m'\cap M}=R$ and $R_n=\overline{\cup_m P_m'\cap M_n}$, 
$n\geq 1$, 
then $(Q\subset R)\subset (N\subset M)$ is a non degenerate commuting square of II$_1$ factors, $Q\subset R\subset^{e_0} R_1 
\subset^{e_1} R_2...$ is its Jones tower and  
$Q'\cap R_n = N'\cap M_n$, $R'\cap R_n=M'\cap M_n$, $\forall n$.

\vskip.05in 

To see this, let $M\vee M^{^{op}}\subset M\underset{e_N}\to{\boxtimes} M^{^{op}}=S$ be the SE inclusion of factors 
associated with $N\subset M$.  By applying recursively 
Lemma 2.8, we obtain a sequence of subfactors $M\supset N = P_0 \supset P_1 ...$ 
such that for each $n\geq 1$ we have 

\vskip.05in 

$(a)$ $P_n\subset P_{n-1}$ is $(P_{n-1}\subset M)$-compatible (thus also $(N\subset M)$-compatible). 

\vskip.05in

$(b)$ $\|E_{(P_n'\cap M)'\cap S}(x)-E_{M'\cap S}(x)\|_2 \leq 2^{-n} \|x\|$,  $\forall x\in P_{n-1}'\cap S$. 

\vskip.05in

Let $Q=\overline{\cup_n P_n'\cap N}$, $R=\overline{\cup_n P_n'\cap M}$. If we denote $S_0=\overline{\cup_n P_n'\cap S}$ , 
then by property $(b)$ above, it follows that $R'\cap S_0=M^{^{op}}$. In particular, $R$ is a II$_1$ factor. By the definitions 
of $Q, R$, it follows that $(Q\subset R)\subset (N\subset M)$ is a commuting square, with $e_0=e_N$ implementing the 
conditional expectation of $R$ onto $Q$ and $Q=\{e_0\}'\cap R$. From $(b)$, one also gets  
$E_{R'\cap S}(e_0)=\lambda 1$. This implies that the algebra  $R^0_1:=\text{\rm sp}R e_0R$ has support $1$ and thus 
any orthonormal basis $\{m_j\}_j$ of $R$ over $Q$  must ``fill up the identity'', i.e.,  $\sum_j m_je_0m_j^*=1$. 
Hence, $(Q\subset R) \subset (N\subset M)$ is in fact a non-degenerate commuting square. Moreover,  ${R^0_1}'\cap S_0=\{e_0\}'\cap M^{^{op}}=
N^{^{op}}$, implying that $R^0_1$ is a II$_1$ factor. Thus, $Q\simeq Qe_0=e_0R^0_1e_0$ is a II$_1$ factor as well, 
and $R^0_n:=vN(R, e_0, ..., e_{n-1})$, $n\geq 1$, is the Jones tower for $Q\subset R$. 

At the same time, 
if for each $n\geq 1$ we define $R_n=\overline{\cup_m P_m'\cap M_n}$, then both this sequence and the sequence 
$Q\subset R \subset R^0_1 \subset ...$  
make (non-degenerate) commuting squares with $N\subset M \subset M_1 \subset ...$, with $R_n^0\subset R_n$. This forces 
$R^0_n=R_n$ and  $R_n$, $N_{n-1}^{^{op}}$ 
be each other's commutant in  $S_0$,  $\forall n \geq 1$. Note that this also implies that for the downward continuation of these towers we have  
$Q'\cap S_0=M_1^{op}$. 

So for the higher relative commutants, we have the equalities  
$$
Q'\cap R_n = (Q'\cap S_0)\cap R_n = (M_1^{^{op}}\cap 
M_n)\cap R_n
$$
$$
=(N'\cap M_n)\cap R_n \subset Q' \cap R_n=M_1^{^{op}} \cap M_n \cap (N_{n-1}^{^{op}})'=N'\cap M_n, 
$$
finishing the proof of Fact 1. 

\vskip.05in 

{\it Fact} 2. Assume $... \subset P_1 \subset P_0 \subset N \subset M$ are as in {\it Fact} 1. If $u_n\in \Cal U(P_n), n\geq 0$, 
and we define $P_n^n=u_0...u_nP_nu_n^*...u_0^*$, then $P_n^n$ is an $(N\subset M)$-compatible subfactor in $P_{n-1}^{n-1}$ and 
$... P_n^n\subset P_{n-1}^{n-1}\subset ... P^1_1\subset P_0^0\subset N\subset M$ is a sequence of factors 
still satisfying the conditions in the statement of {\it Fact} 1. 

\vskip.05in 

 Indeed, for each $k$ the systems of commuting squares of algebras  
 $\{P_n'\cap N \subset  P_n'\cap M \subset ...P_n'\cap M_k\}_n$, 
 $\{{P_n^n}'\cap N\subset {P_n^n}'\cap M \subset ... {P_n^n}'\cap M_k\}_n$   
 (standard $\lambda$-lattices in the sense of [P94a]) are isomorphic  
 via the map $\Phi (x) = \lim_n u_0...u_n x u_n^*.... u_0^*$, $x\in \cup_n P_n'\cap M_k$. 
 Thus,  $\Phi$ implements an isomorphism between $R'\cap R_i \subset  Q'\cap R_i$ 
 and ${R^0}'\cap R_i^0 \subset {Q^0}'\cap R_i^0$, where $R^0_i=\overline{\cup_n {P_n^n}'\cap M_i}$, $i\leq n$. 
 Since the isomorphism $\Phi$ leaves $N'\cap M_i=Q'\cap R_i$ and $M'\cap M_i=R'\cap R_i$ fixed, 
 by equality of dimensions via $\Phi$ it follows that $N'\cap M_i=Q'\cap R_i$, $M'\cap M_i=R'\cap R_i$, $\forall i$.

\vskip.05in 

{\it Fact} 3. Assume $... \subset P_1 \subset P_0 \subset N \subset M$ are as in {\it Fact} 1. 
Then there exist integers $k_0=0 < k_1 < ... $ and unitaries $v_n\in \Cal U(P_{k_{n-1}}), n\geq 0$, 
such that if we define $P_{k_n}^n=v_1...v_nP_{k_n}v_n^*...v_1^*$ 
and let $Q= \overline{\cup_n {P^n_{k_n}}'\cap N}\subset \overline{\cup_n {P^n_{k_n}}'\cap M}=R$, then $Q'\cap M_m=Q'\cap R_m$, 
$R'\cap M_m=R'\cap R_m$, $\forall m$. 

\vskip.05in

To show this, let $\{b_k\}_k\subset (\cup_n M_n)_1$ be a $\| \  \|_2$-dense sequence. We choose recursively $k_m>k_{m-1}, v_m \in \Cal U(P_{k_{m-1}})$ such that 

$$
\|E_{({P^m_{k_m}}'\cap P^{m-1}_{k_{m-1}})'\cap M_m}(b_j)-E_{{P^{m-1}_{k_{m-1}}}'\cap M_m}(b_j)\|_2 \leq 2^{-m}, \forall j\leq m. \tag F3
$$
Assume we made this construction up to $m=n$. Due to Lemma 2.9, (Theorem 0.1 $(a)$ in [P13]) implies that 
if $A_0\subset P_{k_n}^n$ is a finite dimensional abelian von Neumann subalgebra 
with all minimal projections of sufficiently small trace, then there exists $u\in \Cal U(P_{k_n}^n)$ 
such that $\|E_{uA_0u'\cap M_{n+1}}(b_j)-E_{{P_{k_n}^n}'\cap M_{n+1}}(b_j)\|_2 < 2^{-n-1}$, $\forall j\leq n+1$. 
For each $m\geq k_n$ denote $P^n_m = v_0 ... v_n P_m v_n^*... v_0^*$. Since 
$B=\overline{\cup_j {P_j^n}'\cap P_{k_n}^n}$ is diffuse (because it contains a Jones sequence of $\lambda=[M:N]^{-1}$ projections, which generate 
a copy of the hyperfinite II$_1$ factor by [J83]), we may assume $A_0\subset B$, and hence  
$$
\|E_{uBu'\cap M_{n+1}}(x)-E_{{P_{k_n}^n}'\cap M_{n+1}}(b_j)\|_2$$
$$
 < \|E_{uA_0u'\cap M_{n+1}}(b_j)-E_{{P_{k_n}^n}'\cap M_{n+1}}(b_j)\|_2< 2^{-n-1}. 
$$
Since $B$ is a ``limit'' of ${P_j^n}'\cap P_{k_n}^n$, for $j$ sufficiently large we'll still have     
$$
\|E_{({uP^n_ju^*}'\cap P^n_{k_n})'\cap M_{n+1}}(b_j)-E_{{P^n_{k_n}}'\cap M_{n+1}}(b_j)\|_2 < 2^{-n-1}, \forall j\leq n+1.
$$
We choose such a large $j$ and put $k_{n+1}=j$. Letting $v_{n+1} =  v_n^* ... v_1^* u v_1... v_n \in P_{k_n}^n$, 
$P_{k_{n+1}}^{n+1}=v_1...v_{n+1}P_{n+1}v_{n+1}^*...v_1^*=uP^n_{k_{n+1}}u^*$, we see that (F3) is satisfied for $m=n+1$.    

If we now define $Q= \overline{\cup_n {P^n_{k_n}}'\cap N}\subset \overline{\cup_n {P^n_{k_n}}'\cap M}=R$, then 
condition (F3) implies $Q'\cap M_n \subset Q'\cap R_n$, $\forall n$, while Fact 2 implies we have 
$Q'\cap R_n=N'\cap M_n$. The calculations for the relative 
commutants of $R$ are similar, thus finishing the proof. 
\hfill $\square$

\vskip.05in 

Note that the case ``$\Gamma$ finitely generated'' of Theorem 2.1 already follows from Theorem 2.10 above, 
due to the observation we made just before stating Theorem 2.1. But deriving from this the case ``$\Gamma$ infinitely generated'' 
is problematic, as applying it for ``larger and larger'' finitely generated subgroups of $\Gamma$ would involve in the limit 
multiplying infinitely many perturbing unitaries. To deal with this problem we'll use a similar trick, 
but with a ``diagonal embedding'' of $N\simeq N^\sigma$ into the algebra of matrices ``of size $\Gamma$'' over $N$,  $M^\sigma \simeq N \overline{\otimes} \Cal B(\ell^2\Gamma)$.

\proclaim{Lemma 2.11} Let $\{\sigma_g\}_{g\in \Gamma}\subset \text{\rm Aut}(N)$ be an outer action of a group $\Gamma$ on a $\text{\rm II}_1$ factor $N$, 
with $\sigma_e=id_N$, and $v_{g,h}\in \Cal U(N)$ satisfying $\sigma_g\sigma_h = \text{\rm Ad}(v_{g,h})\sigma_{gh}$, 
$v_{e,h}=v_{e,g}=1$, $\forall g,h\in \Gamma$.  
Define $M^\sigma=N \overline{\otimes} \Cal B(\ell^2\Gamma)$ and $N^\sigma=\{\sum_{g} \sigma_g(x)e^0_{gg} \mid 
x\in N\}\subset M^\sigma$, where $\{e^0_{gg'}\}_{g,g'\in \Gamma}\subset \Cal B(\ell^2\Gamma)$ are the usual matrix units.  Let $\Cal R \subset M^\sigma$ 
be a subfactor, with the property that $\{e^0_{gg}\}_g\subset \Cal R$ 
and let $Q\subset \Cal R\cap N^\sigma$ 
be a common subfactor such that $Qe^0_{gg}=e^0_{gg}\Cal Re^0_{gg}$, $\forall g\in \Gamma$.  

\vskip.05in
$1^\circ$ Let $Q_0\subset N$ be the unique subfactor with $Q_0e_{ee}=Qe_{ee}$. 
There exist $w_g\in \Cal U(N)$, $w_e=1$, such that $\sigma'_g:=\text{\rm Ad}(w_g)\circ \sigma_g$ 
satisfy $\sigma'_g(Q_0)=Q_0$, $\forall g\in \Gamma$. 

\vskip.05in
$2^\circ$ If in addition $Q '\cap \Cal R={N^{\sigma}}' \cap \Cal R=\{e^0_{gg}\}_g''$, then ${\sigma'_g}_{|Q_0}$ 
is outer, $\forall g\in \Gamma, g\neq e$, and $v'_{g,h}=w_g \sigma_g(w_h)v_{g,h}w_{gh}^*$ normalize $Q_0$, $\forall g,h$. 

\vskip.05in
$3^\circ$ Let $e\in F_i=F_i^{-1} \nearrow \Gamma$ be the net of finite symmetric subsets,  and 
denote $q^0_i= \sum_{g\in F_i} e^0_{gg}$. Then  
the inclusions of $\text{\rm II}_1$ factors $(Qq^0_i \subset q^0_i\Cal Rq^0_i) \subset (N^\sigma q^0_i \subset q^0_i M^\sigma q^0_i)$ 
make a non-degenerate commuting square $($with respect to the trace preserving conditional expectations$)$, $\forall i$. 

\vskip.05in
$4^\circ$ Denote $M^\sigma_1:=M^\sigma \overline{\otimes} \Cal B(\ell^2\Gamma)\simeq M\overline{\otimes} \Cal B(\ell^2\Gamma)\overline{\otimes} \Cal B(\ell^2\Gamma)$ 
and consider the embedding $j_1:M^\sigma \hookrightarrow M^\sigma_1$ given by $j_1(x) = \sum_g^{-1} \sigma_g(x)e^1_{gg}$, where $\sigma=\sigma \otimes 1$ 
and $\{e^1_{gg'}\}_{g,g'}$ are the matrix units of the $2$nd copy of $\Cal B(\ell^2\Gamma)$. 
If $q_i^1=\sum_{g\in F_i}e^1_{gg}$, then $N^\sigma \subset q^0_iM^\sigma q^0_i \simeq j_1(q^0_iM^\sigma q^0_i)q^1_i \subset q^0_iq^1_i M^\sigma_1q^0_iq^1_i$ 
is a basic construction.  

\vskip.05in
$5^\circ$   With the above notations let $\Cal R^i_1\subset q^0_iq^1_i M^\sigma_1q^0_iq^1_i$ 
be so that $Qq_i^0 \subset q_i^0\Cal Rq_i^0 \hookrightarrow^{j_1} \Cal R^i_1$ is a basic construction. If  
$(Qq^0_i)' \cap \Cal R_1^i=(N^\sigma q^0_iq_i^1)' \cap q_i^0q^1_iM^\sigma_1q_i^0q^1_i$, $\forall i$, then 
$\{{\sigma'_g}_{|Q_0}\}_g$ is an outer action of $\Gamma$ on $Q_0$. If in addition 
$Q'\cap N^\sigma=\Bbb C$ $($equivalently, $Q_0'\cap N=\Bbb C)$, then 
$v'_{g,h}\in Q_0$, $\forall g,h\in \Gamma$, and $\{{\sigma'_g}_{|Q_0}\}_g\subset \text{\rm Out}(Q_0)$ has same $H^3(\Gamma)$-obstruction as $\sigma$. 

\endproclaim
\noindent
{\it Proof.}  $1^\circ$ Identifying $N = N \otimes 1$, let $j_0: N \simeq N^\sigma$  
be the isomorphism defined by the property $j_0(x)e^0_{ee} = xe^0_{ee}, x \in N$. Thus,  $Q_0=j_0^{-1}(Q)$.

Since $\Cal R$ is a (necessarily semifinite) factor and $e^0_{gg}$ belong to $\Cal R$ and have same trace, 
there exist partial isometries $e'_{eg}\in \Cal R$ with left support $e^0_{ee}$, right support $e^0_{gg}$, $\forall g$, and $e'_{ee}=e^0_{ee}$. 
From the way $M^\sigma$ is defined, it follows that there exist unitary elements $w_g\in N$ such that $e_{eg}'=w_ge^0_{eg}$, $\forall g$, with $w_e=1$. 
Let $\sigma'_g=\text{\rm Ad}(w_g)\circ \sigma_g\in \text{\rm Aut}(N)$, $g\in \Gamma$.  Note that for any $x\in N$ 
we have $\sigma'_g(x)e^0_{ee}=e'_{eg}j(x){e'_{eg}}^*$. If in addition $x\in Q_0$, then 
$e'_{eg} j(x){e'_{eg}}^*$ lies in  $e^0_{ee}\Cal Re^0_{ee}$ which is equal to $Qe^0_{ee}=Q_0e^0_{ee}$. Thus, 
$\sigma'_g(Q_0)=Q_0$. 

We clearly have $\sigma'_g\sigma'_h = \text{\rm Ad}(v'_{g,h}) \sigma'_{gh}$, $\forall g, h\in \Gamma$. Since $\sigma'_g$ normalize $Q_0$, it follows that  
$v'_{g,h}$ normalize $Q_0$, $\forall g,h$.
   
$2^\circ$ Note that if $g, g'\in \Gamma$, then $e^0_{gg}({N^\sigma}'\cap M^\sigma)e^0_{g'g'}=0$ 
(resp. $e^0_{gg}(Q'\cap \Cal R)e^0_{g'g'}=0$) is equivalent to $\sigma'_g \neq \sigma'_{g'}$ in $\text{\rm Out}(N)$ (resp. ${\sigma'_g}_{|Q_0} 
\neq {\sigma'_{g'}}_{|Q_0}$ in $\text{\rm Out}(Q_0)$). Applying this to $g'=e$, proves $2^\circ$. 

$3^\circ$ The fact that the inclusions make a commuting square is obvious from $1^\circ$ above. 

$4^\circ$ This is in (Section 5.1.5 in [P91]). 

$5^\circ$ Note that $(Qq^0_iq_i^1)' \cap \Cal R_1^i=(N^\sigma q^0_iq_i^1)' \cap  q_i^0q^1_iM^\sigma_1q_i^0q^1_i$, $\forall i$,  
implies $(Qq^0_i)' \cap q^0_i \Cal Rq^0_i=(N^\sigma q^0_i)' \cap q^0_i M^\sigma q^0_i$, $\forall i$, which implies 
$Q'\cap \Cal R={N^\sigma}'\cap M^\sigma$. So by $2^\circ$, we already know that  ${\sigma'_g}_{|Q_0}$ is outer, $\forall g\neq e$.

Fix $g, h\in \Gamma$ and let $F_i$ be sufficiently large so that $g, h, gh\in F_i$. By taking into account the 
form of the basic constructions $N^\sigma q^0_i \subset q^0_i M^\sigma q^0_i \hookrightarrow q^0_iq^1_0M^\sigma q^0_iq^1_0$ and 
$Qq^0_i \subset q^0_i\Cal R q^0_i\hookrightarrow  \Cal R^i_1$,  as well as 
(Section 5.1.5 in [P97]),  we see like in the proof of $2^\circ$ above that the equality of relative commutants multiplied by $e^0_{ee}e^1_{gh}$ 
from the left and  by $e^0_{gg}e^1_{hh}$from the right implies that $\sigma'_g\sigma'_h$ not equal to  
$\sigma'_{gh}$ in $\text{\rm Out}(N)$ iff ${\sigma'_g}_{|Q_0}{\sigma'_h}_{|Q_0}$ not equal to  
${\sigma'_{gh}}_{|Q_0}$ in $\text{\rm Out}(Q_0)$. This shows that Ad$(v'_{g,h})$ implements inner automorphisms on $Q_0$. 
If in addition $Q_0'\cap N=\Bbb C$, then this forces $v'_{g,h}\in Q_0$, $\forall g, h\in \Gamma$. 
But if $\{v'_{g,h}\}_{g,h}\subset Q_0$, then $\sigma'$ (thus $\sigma$) and $\sigma'_{|Q_0}$ have the same scalar $H^3(\Gamma)$-obstruction. 
\hfill $\square$

\vskip.1in
\noindent
{\it Proof of} Theorem 2.1. By Proposition 1.5.2$^\circ$, we may assume $N$ is separable. 
Also, without loss of generality, we may assume $\sigma_e=id_{N}$. 

Like in Lemma 2.11, we denote $\{e^0_{gh}\}_{g,h\in \Gamma} \subset \Cal B(\ell^2\Gamma)$ the canonical matrix units and  
let $M^\sigma:=N\overline{\otimes} \Cal B(\ell^2\Gamma)$,  $j_0:N \hookrightarrow M^\sigma$ given by  $j_0(x)=\sum_g \sigma_g(x)e^0_{gg}$, 
$N^\sigma=j_0(N) \subset M^\sigma$ . 
Choose an increasing sequence of finite sets  $e\in F_n=F_n^{-1} \subset \Gamma$ 
that exhaust $\Gamma$ and denote $q^0_n=\sum_{g\in F_n} e^0_{gg} $. 

As in Lemma 2.11.5$^\circ$, and with the notations established there, 
the embedding of factors $N\simeq^{j_0} N^\sigma \subset M^\sigma$ 
implements (by induction/reduction) a sequence of inclusions of II$_1$ factors 
$N \simeq N^\sigma q^0_n \subset q^0_nM^\sigma q^0_n \simeq N^{|F_n|}$ whose basic construction identifies, via $j_1( \cdot )q_n^1$, 
with $q^0_nq^1_n M_1^\sigma q^0_nq^1_n $.

For each $n$, the subfactor $N^\sigma q^0_n \subset q^0_nM^\sigma q^0_n$ is a  locally trivial extremal inclusion of II$_1$ factors 
with standard graph given by the Cayley graph of the subgroup $\Gamma_n = \langle F_n\rangle \subset \Gamma$ (see 5.1.5 in [P91]),  
which is amenable.  Using this, we apply Theorem 2.10 
to construct recursively a sequence of subfactors of finite index $N\supset P_0 \supset P_1 \supset...$ such that if we denote 
$Q_0=\overline{\cup_n P_n'\cap N}$, $\Cal R=\overline{\cup_n j_0(P_n)'\cap M^\sigma}$, then $\Cal R$ is an irreducible type II 
subfactor in $M^\sigma$ containing the finite projections $\{e^0_{gg}\}_g$, with $Q:=j_0(Q_0)\subset N^\sigma \cap \Cal R$ a II$_1$ subfactor satisfying the condition 
$Q'\cap \Cal R=Q'\cap M^\sigma={N^\sigma}'\cap M^\sigma=\{e^0_{gg} \mid g\in \Gamma\}''$ 
and with $Qe^0_{gg}=e^0_{gg}\Cal R e_{gg}$, $\forall g$. By Lemma 2.11.$5^\circ$ the inclusion  
$(Qq^0_n \subset q^0_n\Cal R q^0_n)\subset (N^\sigma q^0_n \subset q^0_n M^\sigma q^0_n)$ is a non-degenerate commuting square 
and our construction  will  show that $\Cal R_1^n:=\overline{\cup_m j(P_m)'\cap q^0_nq^1_n M_1^\sigma q^0_nq^1_n}$ is the basic construction 
algebra for the subfactor $(Qq^0_n \subset q^0_n\Cal Rq^0_n)$, $\forall n$, where $j=j_1 \circ j_0$. The desired conclusions will then follow from $2.11.5^\circ$.  

Let $Tr$ denote the semifinite trace  $\tau \otimes Tr_{\Cal B(\ell^2\Gamma)} \otimes Tr_{\Cal B(\ell^2\Gamma)}$ on $M^\sigma_1$ 
and $\|  \ \|_{2,Tr}$ the corresponding $L^2$-norm.  Let $Y=\{y_n\}_n \subset (M^\sigma_1)_1 \cup L^2(M_1^\sigma, Tr)$ 
be a $\| \ \|_{2,Tr}$-dense sequence and denote by $Y_n=\{j_1(q_j^0)q_j^1 y_k j_1(q_j^0)q_j^1 \mid 1\leq j, k \leq n\}$.  We construct the decreasing 
sequence of subfactors $P_m\subset N$ recursively, such that $P_0=N$ and for $m\geq 1$ 

\vskip.05in 

$(i)$ $(j_0(P_m)\subset N^\sigma)q^0_m$ is $(N^\sigma q_m \subset q^0_mM^\sigma q^0_m)$-compatible. 

\vskip.05in 

$(ii)$  $\|E_{j_1((j_0(P_m)'\cap N^\sigma)q^0_m)' \cap M^\sigma_1}(y)-E_{j_1(N^\sigma)'\cap M_1^\sigma}(y))\|_{2,Tr} < 2^{-m}/|F_m|$, $\forall y\in Y_m$. 

\vskip.05in 

If we made this construction up to $m=n$, then by applying Theorem 2.10 to the inclusion of II$_1$ factors 
$j_0(P_n)q^0_{n+1} \subset q^0_{n+1} M^\sigma q^0_{n+1}$ (which has amenable graph by Proposition 2.7, but this is trivial 
here, because this subfactor is in fact locally trivial, with standard graph 
given by a Cayley graph of an amenable subgroup of $\Gamma$, see 5.1.5 in [P97]) we get a  
subfactor $P_{n+1}\subset P_n$ so that its image via $j_0( \ \cdot \ )q^0_{n+1}$ is 
$(N^\sigma q^0_{n+1} \subset q^0_{n+1}M^\sigma q^0_{n+1})$-compatible and such that $(ii)$ above is  
satisfied for $m=n+1$. 

With the sequence $\{P_n\}_n$ this way constructed, define $Q_0$, $Q$ 
and $\Cal R$ as explained above. 
Then $\Cal R$ contains $\{e^0_{gg}\}_g$ and satisfies $Qe^0_{gg}=e^0_{gg}\Cal Re^0_{gg}$, $\forall g$, by construction.

Moreover, condition $(ii)$ shows that 
$$
{j(Q)j_1(q^0_n)q^1_n}'\cap j_1(q^0_n)q^1_nM^\sigma_1j_1(q^0_n)q^1_n
$$ 
$$
=j_1(q^0_n)q^1_n (j_1(N^\sigma)'\cap M^\sigma_1)j_1(q^0_n)q^1_n=j(Q)'\cap \Cal R_1^n.
$$ 
From all this, it also follows that $Q_0'\cap N=\Bbb C$ and that $Q_0$ is AFD, and thus isomorphic to $R$ by ([MvN43]). 

Thus, 2.11.5$^\circ$ applies, to conclude that there exist $\{w_g\}_g \subset \Cal U(N)$ such that 
$\sigma'_g=\text{\rm Ad}(w_g) \sigma_g$ normalizes the irreducible hyperfinite II$_1$ subfactor $Q_0\subset N$, $\forall g$, $v'_{g,h} = w_g\sigma_g(w_h)v_{g,h}w_{gh}^*$ 
belong to $Q_0$, $\forall g,h$, and $\sigma, \sigma'_{|Q_0}$ have the same $H^3(\Gamma)$ obstruction.  
\hfill $\square$

\vskip.05in

We end this section by noticing that the same argument we used in the proof of Theorem 2.1 above shows that any 
``amenable family'' $\Cal G$ of endomorphisms (or Hilbert bimodules) over a II$_\infty$ factor $\Cal N$ normalizes a 
``large'' AFD subfactor $R\overline{\otimes} \Cal B(\ell^2\Bbb N) \subset \Cal N$, 
on which it ``acts faithfully''. To state this result, we need to clarify the terminology and fix some notations. 
Also, the reader should recall that for a properly infinite factor $\Cal N$ (of type II$_\infty$ in our case) one has Connes' well known correspondence  
between a Hilbert $\Cal N$-bimodule $_\Cal N\Cal H_\Cal N$ and a $\Cal N$-endomorphisms $\theta_{\Cal H}$ 
(see [C80]; see also Sec. 1.1 in [P86]). Via this correspondence, the adjoint operation on endomorphisms, $\theta \mapsto \bar{\theta}$,  
is given by $\bar{\theta_{\Cal H}}=\theta_{\bar{\Cal H}}$ and tensor product $\Cal H\overline{\otimes}_{\Cal N}\Cal H'$ 
corresponds to composition $\theta_{\Cal H} \circ \theta_{\Cal H'}$ 
(see e.g., 1.3 in [P86]).  It is this ``nice'' correspondence between Hilbert-bimodules and endomorphisms 
for properly infinite factors that imposes using the framework of II$_\infty$ 
(i.e., infinite amplification of a II$_1$ factor), rather than II$_1$ factors, as algebras on which a category ``acts'' 
(the II$_1$ framework  would instead require considering morphisms between amplifications of the II$_1$ factor, see 1.1 in [P86]).

For us here, if $\Cal N$ is a given II$_\infty$ factor,  a 
{\it concrete C$^*$-tensor category $\Cal G$ of endomorphisms} on $\Cal N$ is a  family of classes (mod 
perturbation by inner automorphisms of $\Cal N$) of unital endomorphisms of $\Cal N$ that contains the $id_\Cal N$ (thus the class of inner automorphisms) 
and satisfies the following properties: $(i)$ it  is closed to the adjoint operation, i.e. if $\theta\in \Cal G$ then  $\bar{\theta} \in \Cal G$; 
$(ii)$ it is closed under composition, i.e, if $\theta, \theta'\in \Cal G $ then $\theta \circ \theta' \in \Cal G$;  
 $(iii)$ Each $\theta\in \Cal G$ dilates the trace $Tr=Tr_\Cal N$ by a finite scalar $1\leq d(\theta)< \infty$, i.e., $Tr\circ \theta = d(\theta) Tr$, 
with the image subfactor  $\theta(\Cal N)\subset \Cal N$ having Jones index  given by the formula $[\Cal N: \theta(\Cal N)]=d(\theta)^2$; $(iv)$ it is closed 
to ``direct sum and subtraction'', in the following sense: $(a)$ if $\theta, \theta'\in \Cal G$ and $v, v'\in \Cal N$ are isometries such that $vv^*+v'{v'}^*=1$ 
then the endomorphism  $\Cal N \ni x\mapsto \theta\oplus \theta'(x):=\text{\rm Ad}(v)\theta(x)+ \text{\rm Ad}(v')\theta'(x) \in \Cal N$ 
belongs to $\Cal G$; $(b)$ if $\theta\in \Cal G$, $0\neq p\in \Cal P(\theta(\Cal N)'\cap \Cal N))$ and $v\in \Cal N$ is an isometry with range $p$, 
then $\text{\rm Ad}(v^*)\theta\in \Cal G$.   

We denote $\text{\rm Irr}(\Cal G)=\{\theta\in \Cal G \mid \theta(\Cal N)'\cap \Cal N=\Bbb C\}$, the family of 
{\it irreducible} elements in $\Cal G$, which we label by the set $K=K_\Cal G$,  as $\{\theta_k\}_{k\in K}$, with $e\in K$ so that $\theta_e=id_\Cal N$ 
and with the adjoint operation implemented by $k\mapsto \bar{k}$ (thus $\theta_{\bar{k}}=\bar{\theta_k}$). 
It is easy to see from the above properties of $\Cal G$ that any $\theta\in \Cal G$ decomposes 
as a finitely supported 
direct sum $\oplus_{k\in K} \ n_k \ \theta_k$ of irreducible endomorphisms in $\Cal G$, with finite multiplicities $n_k \geq 0$ 
(thus $0< \sum_k n_k <\infty$).  

A subset $e\in F_0=\bar{F_0}\subset K$ {\it generates} $\Cal G$, if any $\theta\in \Cal G$ can be obtained from $F_0$ by applying consecutively finitely many 
times the operations $(i), (ii), (iv)$. If $\Cal G$ is generated by a finite such set $F_0\subset K$, one denotes 
by $\Gamma_{\Cal G, F_0}$ (or simply $\Gamma_\Cal G$ when $F_0$ is clear from the context) 
the Cayley-type bipartite graph (or  matrix with non-negative integer entries) $(a_{kk'})_{k, k'\in K}$ with $a_{kk'}$ equal to 
the multiplicity of $\theta_{k'}$ in $(\sum_{i \in F_0} \theta_i) \circ \theta_k$.

We say that such a family $\Cal G$ is {\it amenable} if it satisfies the F\o lner-type condition we mentioned before Proposition 2.7, 
i.e., if given any finite $e\in F_0 = \bar{F_0} \subset K$, and any $\varepsilon >0$, there exists a finite subset $F$ in the  
subcategory $\langle F_0 \rangle \subset \Cal G$ generated by $F_0\subset K$, such that 
$\underset{k'\in \partial_{F_0} F}\to{\sum} d_{k'}^2 < \varepsilon \underset{k\in F}\to{\sum} d_k^2$, where $\partial_{F_0}F$ is the boundary 
of $F$ in $\langle F_0 \rangle$, i.e., the set of all $k'\in K_{\langle F_0 \rangle}\setminus F$ 
such that $a_{k'k}\neq 0$ for some $k\in F$, where $(a_{kk'})_{k,k'}=\Gamma_{\langle F_0 \rangle}$ 
is the graph of $\langle F_0 \rangle$.   By (Theorem 5.3 in [P97]), 
this condition is equivalent to a Kesten-type condition, requiring that $\|\Gamma_{\langle F_0 \rangle}\|^2=
(\sum_{k\in F_0} d_k)^2$ for any finite $e\in F_0=\bar{F_0}\subset K$. 

If $\Cal G$ is a concrete C$^*$-tensor category of endomorphisms on $\Cal N$ 
and $\Cal Q_0\subset \Cal N$ is a II$_\infty$ subfactor, then we say that $\Cal G$ {\it faithfully normalizes} $\Cal Q_0$, if the following conditions are satisfied: 
$(1)$ the trace $Tr$ is semi finite on $\Cal Q_0$; $(2)$ for each endomorphism $\theta\in \Cal G$ there exists $\theta'\in \Cal G$ in 
the same class as $\theta$ such that $\theta'(\Cal Q_0)\subset \Cal Q_0$, and given any other $\theta''\in \Cal G$ in the same class as $\theta$ that satisfies    
$\theta''(\Cal Q_0)\subset \Cal Q_0$, there exists $w_0\in \Cal U(\Cal Q_0)$ with $\text{\rm Ad}(w_0)\circ \theta''_{|\Cal Q_0}=\theta'_{|\Cal Q_0}$; 
$(3)$ If $\theta\in \Cal G$ is so that $\theta(\Cal Q_0)\subset \Cal Q_0$, then $[\Cal Q_0: \theta(\Cal Q_0)]=d(\theta)^2=[\Cal N: \theta(\Cal N)]$ 
and $\theta$ is irreducible non-inner on $\Cal N$ iff $\theta_{|\Cal Q_0}$ irreducible and non-inner on $\Cal Q_0$. 

Note that a  ``concrete C$^*$-tensor category $\Cal G$ of endomorphisms on a factor $\Cal N$'' is an analogue of a ``group 
of outer automorphisms $\Gamma \subset \text{\rm Out}(\Cal N)$''. If one still denotes by $\Cal G$ the underlying 
abstract C$^*$-tensor category, as defined for instance in [NeTu13], then such an object can as well be viewed as an ``outer action of $\Cal G$ 
by endomorphisms on $\Cal N$''. With this interpretation, 
if $\Cal Q_0\subset \Cal N$ is faithfully normalized by $\Cal G$, then the restriction 
map $\theta \mapsto \theta_{|\Cal Q_0}$ is an isomorphism of the underlying abstract rigid C$^*$-tensor categories.  
We say that such a restriction is {\it strongly smooth} if in addition $\theta (\Cal Q_0)'\cap \Cal N = \theta(\Cal N)'\cap \Cal N$, $\forall \theta \in \Cal G$.  

\proclaim{2.12. Theorem} Let $\Cal N$ be a $\text{\rm II}_\infty$ factor and $\Cal G$ an amenable countably generated concrete 
$C^*$-tensor category of endomorphisms  on $\Cal N$, as defined above. Then $\Cal N$ contains 
an AFD $\text{\rm II}_\infty$ subfactor $\Cal R\subset \Cal N$ that's faithfully normalized by $\Cal G$, in the above sense. Moreover, 
if $\Cal N$ is separable, then $\Cal R$ can be chosen so that to satisfy the strong smoothness condition $\theta(\Cal N)'\cap \Cal N 
= \theta(\Cal R)'\cap \Cal N = \theta(\Cal R)'\cap \Cal R$, for all $\theta\in \Cal G$ with $\theta(\Cal R)\subset \Cal R$. 
\endproclaim
\noindent
{\it Proof.} The proof uses the same ideas and follows exactly the same steps as the proof of Theorem 2.1. 
Note first that, as in the proof of Theorem 2.1, it is clearly sufficient to prove the case when $\Cal N$ is separable. 

Next, note that we have a version for endomorphisms of Lemma 2.11.1$^\circ$, as follows:
\vskip.05in

{\it Fact 1.} Let $\Cal N$ be a properly infinite von Neumann factor and $\theta=\{\theta_k\}_{k\in K_0}$ be a  
set of endomorphisms of $\Cal N$, with 
$e\in K_0$ and $\theta_e=id_{\Cal N}$. Define $\Cal M^{\theta}=\Cal N \overline{\otimes} \Cal B(\ell^2K_0)$ and $\Cal N^\theta=\{\sum_{k\in K_0} \theta_k(x)e^0_{kk} \mid 
x\in \Cal N\}\subset \Cal M^\theta$, where $\{e^0_{kk'}\}_{k,k'\in K_0}\subset \Cal B(\ell^2K_0)$ are the usual matrix units.  Let $\Cal R \subset \Cal M^\theta$ 
be a subfactor, with the property that $\{e^0_{kk}\}_k\subset \Cal R$ and $e^0_{kk}$ are mutually equivalent infinite projections in $\Cal R$. Let   
$\Cal Q\subset \Cal R \cap \Cal N^\theta$ 
be a common subfactor such that $\Cal Qe^0_{ee}=e^0_{ee}\Cal Re^0_{ee}$. Then 
there exists a unique subfactor $\Cal Q_0\subset \Cal N$ such that $\Cal Q_0e^0_{ee}=\Cal Qe^0_{ee}$ 
and there exist unitary elements $w_k\in \Cal N$, $k\in K_0$, $w_e=1$, such that $\theta'_k=\text{\rm Ad}(w_k)\circ \theta_k$ 
satisfies $\theta_k'(\Cal Q_0)\subset \Cal Q_0$, $\forall k$. 
\vskip.05in

The proof if the same so we omit it. Let $j_0:\Cal N \simeq \Cal N^\theta \subset \Cal M^\theta$ denote the isomorphism satisfying $j_0(x)e^0_{ee}=we^0_{ee}$, 
$x\in \Cal N$. Thus, $\Cal Q_0=j_0^{-1}(\Cal Q)$. 

In order to make the ideas more transparent, let us first give an argument for the case when $\Cal G$ is generated by a finite subset 
of irreducible endomorphisms $\{\theta_k\}_{k\in K_0}$, with $e\in K_0=\bar{K_0}\subset K$, by exploiting the remark before the statement 
of Theorem 2.1. So let $N$ be a II$_1$ factor such that $N^\infty=\Cal N$, let $M^{\theta, K_0}$ be the amplification of $N$ 
by $\sum_{k\in K_0} d_k$, $\{p_k\}_{k\in K_0}\subset M^{\theta, K_0}$ a partition of $1$ with projections having traces $\tau(p_k)$ 
proportional to $d_k$ and consider its subfactor $N^{\theta,K_0} = \{\sum_{k\in K_0} \theta_k(x)p_k \mid x\in N \simeq p_eM^{\theta,K_0}p_e\}$, 
where $\theta_k$ are endomorphisms chosen in their class so that $\theta_k(p_e)=p_k$, $k\in K_0$. Note that $N^{\theta,K_0}\subset M^{\theta,K_0}$ 
is an extremal inclusion of II$_1$ factors with index $(\sum_{k\in K_0} d_k)^2$ and relative commutant generated by $\{p_k\}_k$. 
We denote $j_0:N=p_e\Cal Np_e \simeq N^{\theta,K_0}$ the identification.  

\vskip.05in

{\it Fact 2.} With the above notation, assume $(Q\subset R)\subset (N^{\theta,K_0}\subset M^{\theta,K_0})$ is a non-degenerate commuting square 
of factors such that $\{p_k\}_k \subset R$ and $p_eRp_e=Qp_e$. Let $(R \subset R_1 \subset ...) \subset (M^{\theta,K_0}\subset M_1 \subset ...)$ 
be the associated tower of commuting squares, obtained by iterating the basic construction. If $Q'\cap R_n = {N^{\theta,K_0}}'\cap M_n$, $\forall n$, 
and one denotes by $Q_0=j_0^{-1}(Q)\subset N$, then $\Cal Q_0=Q_0^\infty \subset N^\infty = \Cal N$ is faithfully normalized by $\Cal G$. 
\vskip.05in

Indeed, by Fact 1 we know that $\theta_k$ can be taken so that $Q=\{\sum_{k\in K_0} \theta_k(x)p_k \mid x\in Q_0\}$ 
and thus so that $\theta_k(\Cal Q_0)\subset \Cal Q_0$. After $\infty$-amplification  
of the tower of commuting squares of factors, the iterated basic construction 
gives rise to consecutive products of endomorphisms of $\Cal N=N^\infty$, $\theta_{k_n} \circ \theta_{k_{n-1}} \circ ... \theta_{k_1}$, for $k_i\in K_0$, 
which all take $\Cal Q_0=Q_0^\infty \subset N^\infty=\Cal N$ into itself. The condition on the compatibility of higher relative commutants amounts to  
$\Cal Q_0$ being faithfully normalized by $\Cal G$.  

Now, if $\Cal G$ is amenable and finitely generated by its irreducible endomorphisms 
indexed by $K_0\subset K$ as above, then $N^{\theta,K_0}\subset M^{\theta,K_0}$ 
has amenable graph, so Theorem 2.10 applies to provide a strongly 
smooth non-degenerate commuting square $(Q\subset R) \subset (N^{\theta,K_0}\subset M^{\theta,K_0})$, with $Q$ and $R$ hyperfinite II$_1$ factors. 
By Fact 2, $Q$ gives rise to a subfactor $\Cal Q_0\simeq Q^\infty$ that's faithfully normalized by $\Cal G$ and satisfies the required strong smoothness 
condition. This proves the statement in the case $\Cal G$ is finitely generated. 

To deal with the infinitely generated  case, we take the 
inclusion $\Cal N^\theta \subset \Cal M^\theta$  in Fact 1 with $K_0=K=K_\Cal G$ and $\{\theta_k\}_{k\in K}=\text{\rm Irr}(\Cal G)$. 
We further embed $\Cal M^\theta = \Cal N \overline{\otimes} \Cal B(\ell^2K)$ into 
$\Cal M_1^\theta:= \Cal N \overline{\otimes} \Cal B(\ell^2K) \overline{\otimes} \Cal B(\ell^2K)$ by $x \mapsto j_1(x)=\sum_k \theta_k(x)e^1_{kk}$, 
where $\{e^1_{kk'}\}_{k,k'\in K}$ are the matrix units in the 2nd copy of $\Cal B(\ell^2K)$. Note that $\Cal N \ni x \in j_1(j_0(x)) \in \Cal M_1^\theta$ 
is given by $j_1(j_0(x))=\sum_{k,k'\in K} \theta_k(\theta_{k'}(x)) e^1_{k'k'}e^0_{kk}$, where $\theta_k=\theta_k \otimes 1$. 

Note that Fact 1 already implies that if $(\Cal Q \subset \Cal R)\subset (\Cal N^\theta \subset \Cal M^\theta)$ is a subinclusion  of factors 
with $\{e^0_{kk}\}_k''={\Cal N^\theta}'\cap \Cal M^\theta \subset \Cal Q'\cap \Cal R$ and $\Cal Qe^0_{ee}=e^0_{ee}\Cal Re^0_{ee}$, 
then there are representants $\theta_k$ such that $\Cal Q_0=j_0^{-1}(\Cal Q)$ satisfies $\theta_k(\Cal Q_0)\subset \Cal Q_0$, $\forall k\in K$. 
Then note the following: 

\vskip.05in

{\it Fact 3}. Assume the above inclusion $(\Cal Q  \subset \Cal R)\subset (\Cal N^\theta \subset \Cal M^\theta)$ is so that:  \newline 
$(a)$ $\Cal Q'\cap \Cal M_1^\theta={\Cal N^\theta}'\cap \Cal M_1^\theta$; $(b)$ there exist finite subsets $e\in K_i=\bar{K_i} \nearrow K$ such 
that if we denote $q^0_i=\sum_{k\in K_i} e^0_{kk}$ then $(\Cal Qq^0_i \subset q^0_i\Cal R q^0_i)\subset (\Cal N^\theta q_i^0 \subset q^0_i\Cal Rq^0_i)$ 
is a non-degenerate commuting square (with respect to the $Tr$-preserving expectations). Then 
$\Cal Q_0=j_0^{-1}(\Cal Q)\subset \Cal N$ is faithfully normalized by $\Cal G$, with its restriction to $\Cal Q_0$ strongly smooth. 

\vskip.05in

The proof of this fact is identical to the proof of $2.11.5^\circ$ above and we leave it as an exercise. 

From this point on, the proof of Theorem 2.12 is similar to the proof of Theorem 2.1. Thus, we let $e\in K_n=\bar{K_n} \nearrow K$ be finite sets, 
denote $D_n=\sum_{k\in K_n} d_k$, let 
$\{p_k\}_{k\in K}\subset \Cal N$ be mutually orthogonal projections of trace $Tr(p_k)=d_k$  such that $1-\sum_{k\in K}$ is an infinite projection. 
By perturbing if necessary each $\theta_k$ by an inner automorphism, we may assume $\theta_k(p_e)=p_k$, $\forall k$. 
We also choose $\theta'_k$ to be in the same class as $\theta_k$ so that 
$\theta'_{k'}(\sum_{k\in K} p_k), k'\in K$, are mutually orthogonal. For each $n$ consider the embedding 
$N=p_e\Cal Np_e  \hookrightarrow^{j^n_0} M^{\theta, K_n}:=N^{D_n}$ by $x \mapsto j^n_0(x)=\sum_{k\in K_n} \theta_k(x)p_k$; as well as the  embedding 
$M^{\theta, K_n} \hookrightarrow^{j^n_1} M_1^{\theta, K_n}:= (M^{\theta, K_n})^{D_n}\simeq N^{D_n^2}$, by $x \mapsto j^n_1(x)=\sum_{k\in K_n} \theta'_k(x)p^n_k$, 
where $p^n_e=\sum_{k\in K_n}p_k$ and $p^n_k = \theta'_k(p^n_e)$. 

Note that these are extremal embeddings of II$_1$ factors with finite index $=D_n^2$ 
and that $j^n_1(j^n_0(N)) \subset j^n_1(M^{\theta,K_n}) \subset M_1^{\theta,K_n}$ is a basic construction,  $\forall n$. 
Moreover, for each $n > m$, the inclusion for $m$ is obtained from the one for $n$ by induction/reduction. We can thus apply the same iterative construction 
of subfactors of finite index $N \supset P_1 \supset P_2 .....$ as in the proof of Theorem 2.1, so that if one denotes $Q_0=\overline{\cup_m P_m'\cap N}$, 
$R^n_0=\overline{\cup_m {j^n_0(P_m)'\cap M^{\theta,n}}}$, $R^n_1= \overline{\cup_m {j^n_1(j^n_0(P_m))' \cap M_1^{\theta,n}}}$, 
then $(j_1^n(j_1^n(Q_0)) \subset j_1^n(R^n_0)\subset R^n_1) \subset (j^n_1(j^n_0(N)) \subset j^n_1(M^{\theta,K_n}) \subset M_1^{\theta,K_n})$ 
is a basic construction of non-degenerate commuting squares that are strongly smooth, $\forall n$. 

But then, it is easy to see that after appropriate amplifications, $Q_0 \subset N$ gives rise to a II$_\infty$ subfactor in $\Cal Q_0 \subset \Cal N$ 
that satisfies the conditions in Fact 3, and is thus faithfully normalized by $\Cal G$, with the strongly smooth condition satisfied. 
\hfill $\square$

\vskip.1in
\noindent
{\bf 2.13. Remarks.} $1^\circ$ The first result about normalizing  a ``large'' hyperfinite subfactor $R$ modulo inner perturbations 
was obtained in [P83], in the case $\Gamma=\Bbb Z$, with $R$ being constructed ``by hand'', 
using Rokhlin towers and an iterative procedure. Shortly after, the question of whether any cocycle action of $\Bbb Z^2$ 
on an arbitrary II$_1$ factor can 
be untwisted was asked in [CJ84]. While we realized at that time that if a similar normalization result could be proved   
for $\Gamma=\Bbb Z^2$ then the problem would reduce to the case $N=R$, where vanishing cohomology 
was just shown in [Oc85],  we could not extend the arguments in [P83] from $\Bbb Z$ to $\Bbb Z^2$, despite much effort. Several years later 
in [P89], we were able to solve this problem  by using tools from subfactor theory. However, the argument in [P89]  could only cover    
groups $\Gamma$ that have a finite set of generators $S\subset \Gamma$ with respect to which $\Gamma$ has 
trivial Poisson boundary (e.g.,  groups with polynomial growth, in particular $\Bbb Z^2$), depending crucially on this condition. 
In retrospect, it is quite surprising that in fact any outer action of any amenable group $\Gamma$ on any II$_1$ factor $N$ normalizes (modulo $\text{\text{\rm Inn}}(N)$) 
an irreducible hyperfinite subfactor, a 
property which turns out to characterize the amenablility of $\Gamma$.

$2^\circ$ Both in the case of outer action $\sigma$ of a group $\Gamma$ on 
$N$ (or more generally of an outer action of a category $\Cal G$ on $\Cal N$) there may be an AFD subalgebra  
$R\subset N$ (resp. $\Cal R\subset \Cal N$) that's normalized by $\Gamma$ (resp. $\Cal G$) but on which the resulting  ``outer action'' 
has either more relations among generators (resulting into an outer action of a quotient of $\Gamma$, resp. $\Cal G$), or fewer relations 
(resulting into an outer action of some $\tilde{\Gamma}$, resp. $\tilde{\Cal G}$, whose quotient is $\Gamma$, resp. $\Cal G$). For instance, 
let $\Bbb F_n \curvearrowright  R$ be a free action of the free group with $n\geq 2$ generators on the hyperfinite II$_1$ factor and $\Bbb F_n \rightarrow \Gamma$ 
a surjective map on a (non-free) group $\Gamma$ with $n$ generators (e.g. $\Gamma = \Bbb Z^n$) and $H$ its kernel. If we let $N=R\rtimes H$ then this gives rise 
to a free cocycle action $(\sigma, v)$ of $\Gamma$ on $N$, which normalizes $R$, but on which it generates a free action of $\Bbb F_n$.  
Similarly, one can take $\Gamma \curvearrowright N=N_0\overline{\otimes} R$ an action that's trivial on $R$ but free on $N$. This normalizes $R$ 
but it implements on it the action of the trivial group.

\heading 3.  Non-vanishing cohomology for amplifications of actions \endheading 

\noindent
{\bf 3.1. Definition} ([P01a]). Let $\Gamma \curvearrowright^\sigma N$ be a free action of a group $\Gamma$ on a type II$_1$ factor $N$. 
Let $p\in N$ be a projection and for each $g\in \Gamma$ choose a partial isometry $w_g\in N$ such that $w_gw_g^*=p$, 
$w_g^*w_g=\sigma_g(p)$ and $w_e=p$. Define $\sigma^p_g\in \text{\rm Aut}(pNp)$ by $\sigma^p_g(x)=w_g \sigma_g(x)w_g^*$, 
$x\in pNp$. Then $\sigma^p$ is a free cocycle action of $\Gamma$ on $pNp$, with $2$-cocycle $v^p_{g,h}=w_g\sigma_g(w_h)w_{gh}^*$, 
$g, h\in \Gamma$.  Moreover, up to cocycle conjugacy $(\sigma^p, v^p, pNp)$ only depends on $\tau(p)=t$, thus 
defining a free cocycle  action $(\sigma^t, v^t)$ of $\Gamma$ on $N^t$, called the {\it amplification of} $\sigma$ {\it by} $t$. 

If in addition $\{w_g\}_g\subset N$ satisfy $w_g\sigma_g(w_h)=w_{gh}$, $\forall g, h\in \Gamma$, 
then $w$ is called a {\it generalized $1$-cocycle} (of support $t$) for $\sigma$, while if this equality holds 
modulo scalars then it is called a {\it weak generalized $1$-cocycle} for $\sigma$ (of support $t$). 

Note that the vanishing (resp. weak-vanishing) of the cocycle $v^t$ amounts to being able to chose the partial isometries 
$\{w_g \mid g\in \Gamma\}\subset N$ so that $w$ be a generalized 1-cocycle (resp. 
weak generalized 1-cocycle) for $\sigma$. 

\proclaim{3.2. Theorem} Let $\Gamma$ be a countable group and $\Gamma \curvearrowright^\sigma R=R_0^{\overline{\otimes} \Gamma}$ the non-commutative Bernoulli $\Gamma$-action with base $R_0\simeq R$.  Let $0<t<1$  
and denote by $\sigma^t$ the free cocycle action obtained by amplifying $\sigma$ by $t$, with 
its $2$-cocycle denoted $v^t$. Assume one of the following properties holds true: 
$(a)$ $\Gamma$ contains an infinite subgroup  with the relative property $(\text{\rm T})$; 
$(b)$ $\Gamma$ contains an infinite subgroup with non-amenable centralizer.  
Then the cocycle $v^t$ is not weak-vanishing. 

\endproclaim 

\noindent
{\it Proof}. Assume by contradiction that the cocycle $v^t$ is weak-vanishing. 

If we are under assumption $(a)$, i.e., if $\Gamma$ has an infinite subgroup with the relative property (T), then (Corollary 4.10 in [P01a]) 
shows that the support $t$ of the generalized weak 1-cocycle $w_g$ must be 1, a contradiction. 

If in turn we are under assumption $(b)$, with $H\subset \Gamma$ an infinite  subgroup with centralizer $H'=\{g\in \Gamma\mid gh = hg, \forall h\in H\}$ 
non amenable, then let $(\alpha, \beta)$ be the s-malleable deformation of $R\overline{\otimes} R$ that commutes with the double action 
$\tilde{\sigma}=\sigma \otimes \sigma: \Gamma \curvearrowright \tilde{R}= R\overline{\otimes} R=(R_0\overline{\otimes} R_0)^{\overline{\otimes} \Gamma}$, as in ([P01a]). 
Denote $M=R\otimes \Gamma$, $M_0=R \rtimes H$, $\tilde{M}= R\overline{\otimes} R \rtimes_{\tilde{\sigma}} \Gamma$
and $\tilde{M}_0=R\overline{\otimes} R \rtimes_{\tilde{\sigma}} H$, with $M$ identified as the 
subfactor $R\otimes 1 \vee \{U_g\}_g$ in $\tilde{M}$, where $\{U_g\}_g\in \tilde{M}$ are the canonical unitaries 
implementing $\tilde{\sigma}$. 

Since $H'$ is non-amenable, the action $\text{\rm Ad}(w_g U_g)$  
of $H'$ on $(p\otimes 1)\tilde{M}_0(p\otimes 1)$ has spectral gap relative to $(p\otimes 1)M_0 (p\otimes 1)$. 
Arguing like in ([P06a]), this implies that $w_hU_h$ and $\alpha_s(w_h)U_h$, $h\in H$, are uniformly $\| \ \|_2$-close, 
for $s\in \Bbb R$ with $|s|$ small. Equivalently,  the map $ \xi \mapsto (w_g\otimes 1) \sigma_g(\xi) \alpha_s(w_g^*\otimes 1)$ gives a unitary representation $\pi_s$ of the group $H$ on the Hilbert space $(p\otimes 1) L^2(R \otimes R)\alpha_s(p\otimes 1)$ which for $|s|$ small has  
the vector $\xi_0=(p\otimes 1)\alpha_s(p\otimes 1)$ almost fixed by $\pi_s(h)$, uniformly in $h\in H$. 

We are thus in exactly the same situation as in the case when $H\subset \Gamma$ is rigid. So the proof of (Corollary 4.10 in [P01a]) applies 
to conclude that the support $t$ of the generalized weak 1-cocycle $\{w_h \mid h\in H\}$ for $\sigma$ must necessarily be equal to $1$, 
a contradiction. 
\hfill $\square$

\proclaim{3.3. Corollary} If $\Gamma$ is a group that contains an infinite subgroup which either has relative property $\text{\rm (T)}$,  
or has non-amenable centralizer in $\Gamma$, then $\Gamma\not\in \Cal V\Cal C_w(R).$
\endproclaim

\proclaim{3.4. Theorem} Let $\Gamma$ be a countable group and $\Gamma \curvearrowright^\rho L(\Bbb F_\Gamma)$ 
the action implemented by the translation from the left by elements $h\in \Gamma$ on the set $\{a_g\}_{g\in \Gamma}$ of 
generators of the free group $\Bbb F_\Gamma$, $\rho_h(a_g)=a_{hg}$, $\forall h, g\in \Gamma$.  
Let $1>t>0$ and denote by $\rho^t$ the free cocycle action obtained by amplifying $\rho$ by $t$, with 
its $2$-cocycle denoted $v^t$. Assume one of the following properties is satisfied: 
$(a)$ $\Gamma$  contains an infinite subgroup with the relative property $(\text{\rm T})$; $(b)$ $\Gamma$  
contains an infinite subgroup with non-amenable centralizer; 
$(c)$ $\Gamma$  contains an infinite amenable group with non-amenable normalizer.  Then the cocycle $v^t$ is not weak-vanishing. 
\endproclaim
\noindent
{\it Proof}. The assumptions $(a)$ and $(b)$ lead to exactly the same argument as in the proof of Theorem 3.2 above, 
by using the ``free s-malleable deformation'' of $L(\Bbb F_\Gamma)$, like in the proof of (Theorem 6.1 in [P01a]). 

So let us assume we are under the assumption $(c)$ and let $H \subset \Gamma$ be an infinite amenable subgroup 
with its normalizer $G\subset \Gamma$ non-amenable.  
Note that $\Bbb F_\Gamma \rtimes \Gamma$ is naturally 
isomorphic to $\Gamma * \Bbb Z$. Denote $N=L(\Bbb F_\Gamma)$, 
$M=N\rtimes \Gamma=L(\Bbb F_\Gamma)\rtimes \Gamma = L(\Gamma * \Bbb Z)$. 
Let $p\in \Cal P(N)$ be a projection of trace $t<1$ and assume $(\rho^p, v^p, pNp)$ has weak vanishing cohomology. 

Thus, if $\{U_g\}_{g\in \Gamma}\subset M$ denote the canonical unitaries implementing $\rho$ on $N=L(\Bbb F_\Gamma)$, 
then there exist partial isometries $w_g\in N$ of left support $p$ such that $U'_g=w_gU_g\subset \Cal U(pNp)$ gives a projective representation 
of $\Gamma$ with scalar $2$-cocycle $\mu$. 
In particular, $\{U'_g \mid g\in G\}''$ gives an embedding of $L_\mu(H)\subset L_\mu(G)$ into $pMp$.

Note that $L_\mu(H)$ is amenable diffuse, $L_\mu(G)$ 
has no amenable direct summand and $L_\mu(H)$ is regular in $L_\mu(G)$. Thus, by (Corollary 1.7 in [I13]) it follows that $L_\mu(G)\prec_M L(\Gamma)$. 

But $L_\mu(G)\prec_M L(\Gamma)$ implies that the ``free'' malleable deformation $(\alpha, \beta)$  of $M\subset \tilde{M}=L(\Bbb F_\Gamma * \Bbb F_\Gamma) \rtimes \Gamma$ is uniform on $\{U'_g \mid g\in G\}$ (because $L(\Gamma)$ is in the fixed point algebra 
of the malleable deformation and $L_\mu(G)$ is subordinated to $L(\Gamma)$). Thus, for each $s\in \Bbb R$, 
the $G$-representation $ \xi \mapsto (w_g * 1) \sigma_g(\xi) \alpha_s(w_g^* * 1)$ gives a unitary representation $\pi_s$ of the group $G$ on the Hilbert space $(p\otimes 1) L^2(N * N)\alpha_s(p\otimes 1)$ which for $|s|$ small has  
the vector $\xi_0=(p * 1)\alpha_s(p * 1)$ almost fixed by $\pi_s(g)$, uniformly in $g\in G$.   As 
in the proof of (Theorem 4.1 in [P01a]) this shows that we must necessarily have $t=1$, a contradiction. 
\hfill $\square$

\vskip.1in 
\noindent
{\bf 3.5. Remark.} Let $(N_0, \tau_0)$ be an arbitrary tracial von Neumann algebra $\neq \Bbb C$ 
and $\Gamma$ a group that satisfies one of the conditions $(a)$ or $(b)$ in Theorem 3.2. 
If $\Gamma \curvearrowright N=N_0^{\overline{\otimes} \Gamma}$ is the Bernoulli action with base $N_0$ 
and $0<t<1$, then by using the s-malleable deformation of $N=N_0^{\overline{\otimes} \Gamma}$ in ([I06]) 
in combination with the argument in the proof of Theorem 3.2, one obtains that the 
$t$-amplification $(\sigma^t, v^t, N^t)$ of $\sigma$ has non-vanishing cohomology. Similarly, 
by using the s-malleable deformation in  [IPeP05] one can obtain a generalization of Theorem 3.4 
for the free-Bernoulli action $\Gamma \curvearrowright^\rho N=N_0^{* \Gamma}$. 

\heading 4. Non-vanishing cohomology for cocycle actions on $L(\Bbb F_\infty)$  \endheading 

\noindent
{\bf 4.1. Definition} ([CJ84]). Let $\Gamma$ be an infinite countable group with a set of generators $S\subset \Gamma$ and assume 
$(\Gamma, S) \neq (\Bbb F_S, S)$. Let $\pi: \Bbb F_S \rightarrow \Gamma$ be the unique group morphism taking 
the free generators $S$ of $\Bbb F_S$ onto $S\subset \Gamma$ (so our assumption is equivalent to $\pi$ not being 1 to 1).  Then  
$\text{\rm ker}\pi \simeq \Bbb F_\infty$ and $\Bbb F_S$ has infinite conjugacy classes relative to $\text{\rm ker}\pi$. 
Thus, the inclusion of II$_1$ factors $L(\text{\rm ker}\pi)=N \subset M=L(\Bbb F_S)$ is irreducible and regular, 
with $\Cal N_M(N)/\Cal U(N)=\Gamma$. Consequently,  $N\subset M$ is a crossed product inclusion of the form 
$L(\Bbb F_\infty) \subset L(\Bbb F_\infty) \rtimes_{(\sigma_\pi, v_\pi)} \Gamma$, for some free cocycle action $(\sigma_\pi, v_\pi)$ 
of $\Gamma$ on $L(\Bbb F_\infty)$, that we'll call the {\it Connes-Jones cocycle action} associated with $\pi$, 
with $v_\pi$ its 2-{\it cocycle}, called the  {\it CJ-cocycle}. 

\vskip.05in

We summarize here some straightforward consequences of results from ([CJ84], [P01a], [O03], [P06b], [OP07]), but which 
are stated in a manner pertaining to our vanishing cohomology problem: 

\proclaim{4.2. Theorem} 
Let $(\Gamma, S)$, $\pi: \Bbb F_S \rightarrow \Gamma$, $\Gamma \curvearrowright^{(\sigma_\pi, v_\pi)} L(\Bbb F_\infty)$ 
be as in $4.1$. Assume $\Gamma$ satisfies one of the following conditions: 
$(a)$ it does not have 
Haagerup property $($e.g., it contains an infinite subset with relative property $\text{\rm (T)})$; $(b)$ it has Cowling-Haagerup 
invariant $\Lambda(\Gamma)$ larger than $1$; $(c)$ it has an infinite subgroup with non-amenable centralizer; $(d)$ it has 
an infinite amenable subgroup with non-amenable normalizer.  
Then the CJ-cocycle $v_\pi$ is not weak-vanishing. Thus, if $\Gamma$ satisfies any of these properties then 
$\Gamma \not\in \Cal V\Cal C_w(L(\Bbb F_\infty))$. 
\endproclaim
\noindent
{\it Proof}. If $v_\pi$ is weak-vanishing, then we can choose 
representatives $\{U_g \mid g\in \Gamma\}$ in $\Cal N_M(N)$ 
such that $U_gU_h=\mu_{g,h}U_{gh}$, $\forall g, h\in \Gamma$, 
for some scalar $2$-cocycle  $\mu\in H^2(\Gamma)$. Thus, $L_\mu(\Gamma) \subset M=L(\Bbb F_S)$. 
Also, the unitaries $U_g^{op}\in L(\Bbb F_S)^{op}\simeq L(\Bbb F_S)$ satisfy $U_gU_h=\overline{\mu_{g,h}}U_{gh}$, 
implying that $\{U_g\otimes U_g^{op}\}_g''\simeq L(\Gamma)$ embeds into $L(\Bbb F_S)\overline{\otimes} L(\Bbb F_S)$. 

Since $\Bbb F_S$ has Haagerup property and $\Lambda(\Bbb F_S)=1$ (cf [H79]), the II$_1$ factor $L(\Bbb F_S)$ has Haagerup property 
and $\Lambda(L(\Bbb F_S))=1$, implying that $L(\Bbb F_S)\overline{\otimes} L(\Bbb F_S)$ and all its 
von Neumann subalgebras have these properties as well (see e.g. [P01b] and [OP07]). In particular, $L(\Gamma)$ has these properties, 
thus $\Gamma$ has Haagerup property and $\Lambda(\Gamma)=1$.

If in turn $\Gamma$ has an infinite subgroup $H$ with non-amenable centralizer $H'$, then 
by (Remark $3^\circ$  in $\S 4$ of [P06b]) $L_\mu(\Gamma)$ would contain a diffuse 
von Neumann subalgebra $B$ with non-amenable centralizer. But then $B\subset M=L(\Bbb F_S)$ 
has non-amenable centralizer, contradicting the solidity of free group factors ([O03]). 

If $\Gamma$ has an infinite amenable subgroup $H \subset \Gamma$ 
with non-amenable normalizer $G\subset \Gamma$, then the normalizer of $B=L_\mu(H)$ in $M=L(\Bbb F_S)$ 
contains all $\{U_g\}_{g\in G}$. Since $\{U_g\}_g''\simeq L_\mu(G)$ is non-amenable, this contradicts  
 the strong solidity of free group factors ([OP07]).  

\hfill $\square$

\vskip.05in 
\noindent
{\bf 4.3. Notations}. Given two discrete groups $\Gamma, \Lambda$, one writes $\Gamma \leq_{_{W^*}} \Lambda$ 
whenever $L(\Gamma)$ can be embedded (tracially) into 
$L(\Lambda)$. Denote $\text{\rm W}^*_{leq}(\Lambda)$ the class of all groups $\Gamma$ 
that can be subordinated this way to $\Lambda$, i.e.,  
$\text{\rm W}^*_{leq}(\Lambda) = \{\Gamma \mid \Gamma \leq_{_{W^*}} \Lambda\}$).

\proclaim{4.4. Corollary} One has $\Cal V\Cal C\subset \Cal V\Cal C(L(\Bbb F_\infty)) \subset \text{\rm W}^*_{leq}(\Bbb F_2)$.  Thus, if $\Gamma \in \Cal V\Cal C$, 
then $\Gamma$ has Haagerup property, Cowling-Haagerup constant equal to $1$, any infinite subgroup of $\Gamma$ has amenable centralizer, and any infinite 
amenable subgroup has amenable normalizer. 
\endproclaim
\noindent

\noindent
{\bf 4.5. Remarks} $1^\circ$ The above criteria for $\Gamma$ not to be in $\text{\rm W}^*_{leq}(\Bbb F_2)$ (and thus not in $\Cal V\Cal C(L(\Bbb F_\infty))\supset \Cal V\Cal C$) 
are in fact not stated in their optimal form. Thus, the results in [OP07] show that in order for 
$\Gamma$ not to be in $\text{\rm W}^*_{leq}(\Bbb F_2)$, it is enough that $L(\Gamma)$ is not {\it strongly solid}, 
i.e., that $L(\Gamma)$ merely has a diffuse amenable von Neumann subalgebra 
whose normalizer in $L(\Gamma)$ generates a non-amenable von Neumann algebra. The list in 4.2 is also not exhaustive. For instance, 
by [Pe07] it follows that 
if $\Gamma \in \text{\rm W}^*_{leq}(\Bbb F_2)$, then $L(\Gamma)$ needs to be $L^2$-{\it rigid}, while a result of Ozawa (see e.g. [BrO08]) shows that $\Gamma$ needs to be {\it exact}. One should also 
note that in the proof of Theorem 4.2 above we showed that an embedding $L_\mu(\Gamma)\hookrightarrow  L(\Bbb F_2)$, for some $\mu\in H^2(\Gamma)$, 
gives rise to an embedding $L(\Gamma) \hookrightarrow L(\Bbb F_2 \times \Bbb F_2)$, by simply doubling the 
canonical unitaries $\{u_g\mid g\in \Gamma\}\subset L_\mu(\Gamma)\subset L(\Bbb F_2)$, i.e., by 
taking $L(\Gamma)\simeq \{u_g \otimes u_g^{op}\mid g\in \Gamma\}'' \subset L(\Bbb F_2)\overline{\otimes} L(\Bbb F_2)=L(\Bbb F_2 \times \Bbb F_2)$. 

\vskip.05in

$2^\circ$ It is reasonable to expect that $\Cal V\Cal C= \Cal V\Cal C(L(\Bbb F_\infty))$, and more specifically that the 
CJ-cocycles are in some sense the ``worse possible'', ie., if any such cocycle vanishes for some group 
$\Gamma$, then $\Gamma\in \Cal V\Cal C$. We also believe that $\Cal V\Cal C=\Cal V\Cal C_w$. 

\vskip.05in

$3^\circ$ Given a group $\Lambda$, denote by $\text{\rm ME}(\Lambda)$ the class 
of groups $\Gamma$ that are measure equivalent (ME)  to $\Lambda$ and by $\text{\rm ME}_{leq} (\Lambda)$  
the class of groups $\Gamma$ that have a free  m.p. action which can be realized as a sub equivalence 
relation of a  free ergodic m.p. $\Lambda$-action. It would be interesting to explore the possible correlations 
between the classes $\Cal V\Cal C$, $\text{\rm W}^*_{leq}(\Bbb F_2)$, $\text{\rm ME}_{leq}(\Bbb F_2)$, etc.  In this respect, one should point out that 
while $\text{\rm W}^*_{leq}(\Bbb F_2)$, $\text{\rm ME}_{leq}(\Bbb F_2)$ are obviously  ``hereditary'' classes (i.e., if $\Gamma$ belongs to 
any of them, then all subgroups of $\Gamma$ belong too), we could not prove such hereditarity for 
$\Cal V\Cal C$ (cf. Remark $1.6.2^\circ$). See also  Section 7 in [PeT10] for more comments on $\text{\rm ME}_{leq}(\Bbb F_2)$ and its relations 
to $\text{\rm W}^*_{leq}(\Bbb F_2)$. 
One should also note that $\text{\rm ME}_{leq}(\Bbb F_2)$ consists of groups $\Gamma$ that are ME to either 
$\Bbb Z=\Bbb F_1$, $\Bbb F_2$, or $\Bbb F_\infty$, i.e., $\text{\rm ME}_{leq}(\Bbb F_2)=\text{\rm ME}(\Bbb Z)\cup \text{\rm ME}(\Bbb F_2)\cup \text{\rm ME}(\Bbb F_\infty)$ 
(cf. [G04], [Hj04]). 

\vskip.05in

$4^\circ$ We do not know of any examples of groups in $\Cal V\Cal C$, $\text{\rm W}^*_{leq}(\Bbb F_2)$ other than 
amalgamated free products of amenable groups over finite groups. The approach in 4.2  indicates that these two classes 
may coincide (perhaps with $\text{\rm ME}_{leq}(\Bbb F_2)$ as well). 

An intriguing class of groups  
that are known from [G04] to belong to $\text{\rm ME}_{leq}(\Bbb F_2)$ (in fact, even to 
$\text{\rm ME}(\Bbb F_2)$), are 
the free products of finitely many copies of $\Bbb F_2$ with amalgamation over the subgroup $\Bbb Z\subset \Bbb F_2$ generated 
by the commutator $aba^{-1}b^{-1}$ (where $a, b$ are the generators of $\Bbb F_2$). Gaboriau 
conjectured that in fact any amalgamated free product of $\Bbb F_{k_i}$, $2\leq k_i \leq \infty$, $i\in I$ (with $I$ finite or $I=\Bbb N$), 
over some $\Bbb Z\hookrightarrow \Bbb F_{k_i}$ which is maximal abelian 
in the corresponding $\Bbb F_{k_i}$, $\forall i$, is in ME$_{leq}(\Bbb F_2)$. 

Thus, according to the above speculations, the groups $\Bbb F_{k_1} *_\Bbb Z \Bbb F_{k_2} *_\Bbb Z ...$, with $\Bbb Z$ maximal abelian 
in each $\Bbb F_{k_i}$, should belong to $\Cal V\Cal C$ as well. 
However, we were not able to prove this 
for any such example, except of course the case when any subgroup $\Bbb Z$ is freely complemented in $\Bbb F_{k_i}$.  
Related to this, we pose here the following 

\vskip.05in 
\noindent{\bf Question:} {\it Let $\Bbb F_n \curvearrowright^\sigma N$ be a free action of a free group of rank $n$ on a $\text{\rm II}_1$ factor $N$ 
and let $\Cal W$ denote the set of all $1$-cocycles $w$ for $\sigma$ $($i.e., maps $w:\Bbb F_n \rightarrow \Cal U(N)$ 
satisfying $w_g\sigma_g(w_h)=w_{gh}$, $\forall g, h\in \Bbb F_n)$. 
Let $\Bbb Z\subset \Bbb F_n$ be a maximal abelian subgroup, generated by some element $g\in \Bbb F_n$. Is it then true that the set 
$\{w_g \mid w\in \Cal W\}$ coincides with the unitary group $\Cal U(N)$}?  
\vskip.05in 

Taking into account the way the 1-cocycles $w\in \Cal W$ are constructed, from $n$-tuples 
of unitaries in $N$ that are taken as perturbations of the canonical unitaries $U_1, ..., U_n \in N\rtimes_\sigma \Bbb F_n$ 
that implement $\sigma_{a_1}, ...., \sigma_{a_n}$ (where $a_1, ..., a_n$ are the generators of $\Bbb F_n$), 
it immediately follows  that this question has an affirmative answer in the case $\Bbb Z$ is freely complemented in $\Bbb F_n$. 
It is also trivial to see that if the statement holds true, then all of the above Gaboriau groups $\Gamma=\Bbb F_{k_1} *_\Bbb Z \Bbb F_{k_2} *_\Bbb Z  ...$ 
lie in $\Cal V\Cal C$. 

It is not known whether these groups are in $\text{\rm W}^*_{leq}(\Bbb F_2)$ either. In fact, deciding that a  group 
$\Gamma$ satisfies $L(\Gamma) \hookrightarrow L(\Bbb F_2)$ is at least as interesting as  deciding that it has property $\Cal V\Cal C$. 
So the fact that $\Cal V\Cal C \subset \text{\rm W}^*_{leq}(\Bbb F_2)$ gives 
another strong motivation for  proving  the  universal vanishing cohomology property for various groups. 

The above question can also be stated for free measure preserving actions on the probability measure space, $\Bbb F_n \curvearrowright (X,\mu)$, 
by simply replacing $N$ by $A=L^\infty(X,\mu)$ throughout that statement. Besides answering this question, it would be interesting to know if an affirmative 
answer  would imply Gaboriau's conjecture that the groups $\Bbb F_{k_1} *_\Bbb Z \Bbb F_{k_2} *_\Bbb Z  ...$ belong to ME$_{leq}(\Bbb F_2)$.

\vskip.05in

$5^\circ$ It has been conjectured by Peterson and Thom (see end of  Sec. 7 in [PeT10]) that if two amenable von Neumann subalgebras 
$B_1, B_2$ of the free group factor $L(\Bbb F_2)$ have diffuse intersection, then $B_1\vee B_2$ should follow amenable. 
There has been an accumulation of evidence towards this fact being true (e.g., [P81], [Ju06], [Pe07]). For us here, this would imply that 
if $\Gamma \in \text{\rm W}^*_{leq}(\Bbb F_2)$ is generated by amenable subgroups $\Gamma_1, \Gamma_2\subset \Gamma$ with $H=\Gamma_1 \cap \Gamma_2$ infinite, then $\Gamma$ must be amenable. Thus, if $\Gamma = \Gamma_1 *_H \Gamma_2$ then $\Gamma \not\in \text{\rm W}^*_{leq}(\Bbb F_2)$, 
unless either $H$ is finite, or $[\Gamma_1:H]\leq 2, [\Gamma_2:H]\leq 2$. In particular, $\Cal V\Cal C(L(\Bbb F_\infty))$ should not contain 
such groups either. 

\vskip.05in

$6^\circ$ We expect that $\Cal V\Cal C(R)$ is equal to $\Cal V\Cal C(L(\Bbb F_\infty))$. This fact suggests 
various new statements in deformation-rigidity for factors arising from Bernoulli actions. For instance, it should be possible to 
prove that if a group $\Gamma$ does not have  Haagerup property, or if it has an infinite amenable subgroup with non-amenable normalizer, 
then  some of the W$^*$-rigidity results in [P01a], [P03], [P06a] should be true. Combining this conjecture  with the Peterson-Thom conjecture and remark $5^\circ$ above, 
this also suggests that $\Cal V\Cal C(R)$ doesn't contain any non-amenable group $\Gamma$ that can be generated by amenable subgroups $\Gamma_1, \Gamma_2$ 
with $H=\Gamma_1 \cap \Gamma_2$ infinite (in particular 
$\Gamma=\Gamma_1 *_H \Gamma_2$). But the obstruction in this case should be of a completely different nature. 
The II$_1$ factors arising from Bernoulli actions of such groups (and more generally  from non-amenable groups 
$\Gamma$ generated by $n\geq 2$ amenable groups with infinite intersection) 
may actually have additional W$^*$-rigidity properties, providing a new class 
of factors on which deformation-rigidity techniques should be tested, in the spirit of ([P01a], [P03], [P06a], [IPeP05], [PV12]).  
 
\vskip.05in 

$7^\circ$ Note that  the $W^*$-algebra version of von Neumann's conjecture on whether any 
non-amenable group $\Lambda$ contains a copy of $\Bbb F_2$, amounts to whether for any non-amenable $\Lambda$ one has 
$\Bbb F_2 \leq_{_{\text{\rm W}^*}} \Lambda$. Note also that by   
the Gaboriau-Lyons result in [GL07] 
one indeed has $\Bbb F_2 \leq_{_{\text{\rm W}^*}} \Bbb Z \wr \Lambda$ for any non-amenable $\Lambda$, while by [OP07]  
it follows that if $\Bbb Z \wr \Lambda \leq_{_{\text{\rm W}^*}} \Bbb F_2$ then $\Lambda$ must be amenable.

\heading 5. A related characterization of amenability   \endheading 

We prove in this section that the ``normalization'' property for cocycle $\Gamma$-actions in Theorem 2.1 
can only be true when the group $\Gamma$ is amenable. More precisely, for any non-amenable $\Gamma$  
we exhibit examples of embeddings $\Gamma \subset \text{\rm Out}(L(\Bbb F_\infty))$ 
which admit no lifting to $\text{\rm Aut}(L(\Bbb F_\infty))$ that normalizes a hyperfinite subfactor of $L(\Bbb F_\infty)$. 

\proclaim{5.1. Theorem} Let $\Gamma$ be a countable group and  $\Gamma \curvearrowright^\sigma L(\Bbb F_\Gamma)$ 
the action implemented by the translation from the left by elements $h\in \Gamma$ on the set $\{a_g\}_{g\in \Gamma}$ of 
generators of the free group $\Bbb F_\Gamma$. Let $M=L(\Bbb F_\Gamma)\rtimes \Gamma$ with $\{U_g\mid g\in \Gamma\}\subset M$ 
the canonical unitaries implementing $\sigma$. The following conditions are equivalent: 

\vskip.05in 
$(a)$ $\Gamma$ is amenable. 

$(b)$ $L(\Bbb F_\Gamma)$ contains a hyperfinite subfactor $R$ with $R'\cap M=\Bbb C$ 
for which there exist $\{w_g\}_{g\in \Gamma} \subset \Cal U(L(\Bbb F_\Gamma))$ 
with the property that $U'_g=w_gU_g$, $g\in \Gamma$, 
normalize $R$, implement a free $\Gamma$-action on it and satisfy $U'_gU'_h=U'_{gh}$, $\forall g, h\in \Gamma$.

$(c)$  $L(\Bbb F_\Gamma)$ contains a diffuse AFD von Neumann subalgebra $B$ 
for which there exist $\{w_g\}_{g\in \Gamma} \subset \Cal U(L(\Bbb F_\Gamma))$ 
with the property that $\text{\rm Ad}(w_g)\circ \sigma_g$ normalize $B$, $\forall g$. 
\endproclaim
\noindent
{\it Proof}. By Theorem 2.2 we have $(a) \Rightarrow (b)$, while $(b) \Rightarrow (c)$ is trivial.  

To see that $(c)\Rightarrow (a)$, note first that one has a natural identification between $M=L(\Bbb F_\Gamma) \rtimes_\sigma \Gamma$ 
and $L(\Gamma) * L(\Bbb Z)$. Then note that  by [OP07], the von Neumann algebra $B_1$ generated by the normalizer of $B$ in $L(\Bbb F_\Gamma)$ 
then $B_1$ is amenable, thus AFD by [C75]. As in part $(b)$, denote $U'_g = w_g U_g\in M$, $g\in \Gamma$. Note that $\{U'_g\}_g$  
normalize $B_1$ and implement a cocycle action of $\Gamma$ on $B_1$, i.e., one has $U'_{gh}(U_g'U'_h)^{-1}\in B_1$, $\forall g, h\in \Gamma$. 

We first show that $\Gamma$ non-amenable 
implies  $P=B_1 \vee \{U'_g \mid g\in \Gamma\}$ is non-amenable. We will prove this by contradiction, by showing that if $P$ is amenable, 
then $\Gamma$ follows amenable. 

To see this, let $P \subset \langle P, e_{B_1} \rangle$ be the basic construction for $B_1\subset P$ and note that $f_g:=U'_g e_B {U'_g}^*\in 
\langle P, e_B \rangle$, $g\in \Gamma$, are mutually orthogonal projections summing up to $1$ and generating an atomic abelian 
von Neumann subalgebra $\Cal A \subset \langle P, e_{B_1} \rangle$ naturally isomorphic to $\ell^\infty\Gamma$. Moreover,  $\{\text{\rm Ad}(U'_g)\}_g$ 
normalizes $\Cal A\simeq \ell^\infty\Gamma$ implementing on it 
the $\Gamma$-action by left translation. If $P$ is amenable then there exists a state $\varphi$ on 
$\Cal B(L^2P)$ that has $P$ in its centralizer. Then $\psi=\varphi_{|\Cal A}$ is a state on $\Cal A\simeq \ell^\infty\Gamma$ 
that's invariant to $\{\text{\rm Ad}U'_g\}_g$, i.e., to left translations by $\Gamma$, showing that $\Gamma$ is amenable. 

Now, if $P\subset M\simeq L(\Gamma) * L(\Bbb Z)$ is non-amenable, 
then by (Corollary 1.7 in [I13]) it follows that $P\prec_M L(\Gamma)$, in particular $B\prec_M L(\Gamma)$. But by applying (Corollary 2.3 in [P03]), 
it is trivial to see that $L(\Bbb F_\Gamma)$ (an algebra on which $\{U_g\}_{g\in \Gamma}$ acts) 
has no diffuse von Neumann subalgebra that can be subordinated $\prec_M$ to $L(\Gamma)=\{U_g\}_g''$, a contradiction.  
\hfill $\square$

\proclaim{5.2. Theorem} Let $\Gamma$ be a countable group with a set of generators $S\subset \Gamma$ and  
corresponding surjective morphism $\pi: \Bbb F_S\rightarrow \Gamma$, with kernel $\text{\rm ker}\pi\simeq \Bbb F_\infty$. Let 
$L(\Bbb F_\infty) = N \subset M=L(\Bbb F_S)$ be the associated irreducible, regular  
inclusion of free group factors, which satisfies $\Cal N_M(N)/\Cal U(N)\simeq \Gamma$, and denote 
by $\Gamma \curvearrowright^{\sigma_\pi} N$ the corresponding 
cocycle $\Gamma$-action.   The following conditions are equivalent: 

\vskip.05in 
$(a)$ $\Gamma$ is amenable. 

$(b)$ $N$ contains a hyperfinite subfactor $R \subset N$ with $R'\cap M=\Bbb C$ 
for which there exist $\{w_g\}_g \subset \Cal U(N)$ 
with the property that $U'_g=w_gU_g$, $g\in \Gamma$, 
normalize $R$, implement a free action on it, and satisfy $U'_gU'_h=U'_{gh}$, $\forall g, h\in \Gamma$.

$(c)$ $N$ contains a diffuse AFD von Neumann 
subalgebra $B$ such that $\forall g\in \Gamma$, 
$\exists w_g\in \Cal U(N)$ 
with the property that $\text{\rm Ad}w_g \circ \sigma_g$ normalizes $B$.

\endproclaim
\noindent
{\it Proof}. Theorem 2.2 shows that $(a)\Rightarrow (b)$ and $(b)\Rightarrow (c)$ is trivial. If 
$(c)$ holds but we assume $\Gamma$ is non-amenable, then by (Theorem 3.2.4 in [P86]; see also 
the direct argument in the proof of Theorem 5.2 above) the von Neumann algebra generated 
by $B$ and its normalizer in $M=L(\Bbb F_S)$ is non-amenable, contradicting the strong solidity of 
the free group factors ([OP07]).  
\hfill $\square$

\vskip.05in 
\noindent
{\bf 5.3. Remark.} The dichotomy amenable/non-amenable in the above results can probably be extended to 
cover the converse to Theorem 2.10 as well. Thus, it should be true that if $\Cal G$ is a non-amenable 
standard $\lambda$-lattice, then there exists an extremal inclusion of separable II$_1$ factors $N\subset M$ with standard invariant equal to 
$\Cal G$ in which one cannot embed with non-degenerate strongly smooth commuting squares any inclusion 
of hyperfinite II$_1$ factors $Q\subset R$  (as before, {\it strongly smooth} commuting square 
inclusion $(Q\subset R) \subset (N\subset M)$ means that it is non-degenerate 
and satisfies $Q'\cap R_n=Q'\cap M_n = N'\cap M_n$, $R'\cap R_n = R'\cap M_n = M'\cap M_n$, $\forall n$). 

Similarly, it should be true that if $\Cal G$ is a non-amenable countable rigid $C^*$-tensor category, then $\Cal G$ admits 
an action (by endomorhisms) on a II$_\infty$ factor $M^\infty$ which has no hyperfinite II$_\infty$ subfactors $R^\infty$ with normal 
expectation that are left invariant by $\Cal G$ (modulo inner perturbations) and on which $\Cal G$ acts outerly as a $C^*$-tensor category. 

In both statements, the obvious candidate for a proof is the canonical inclusion  $N^{\Cal G}(L(\Bbb F_\infty)) \subset M^{\Cal G}(L(\Bbb F_\infty))$, 
from ([P94a]), which has $L(\Bbb F_\infty)$ as ``initial data'' and $\Cal G$ as standard invariant. Note that the resulting factors 
$N=N^{\Cal G}(L(\Bbb F_\infty)), M=M^{\Cal G}(L(\Bbb F_\infty))$ were in fact shown to be isomorphic to $L(\Bbb F_\infty)$ in ([PS01]). 
If one assumes by contradiction that there does exist a hyperfinite inclusion $Q\subset R$ with $\Cal G$ as 
standard invariant and which can be embedded with strongly smooth commuting square into $N\subset M$ 
and one takes the associated SE inclusions, then deformation/rigidity  arguments in the style of the proofs of 5.1 and 5.2 above should 
contradict the non-amenability of $\Cal G$. A study case is when $\Cal G$ is the Temperley-Lieb-Jones $\lambda$-lattice 
$\Cal G_\lambda$ of index 
$\lambda^{-1} >4$. One difficulty in proving such a result is that so far (relative) strong solidity  results can  only 
say something about  {\it normalizers}  
of diffuse AFD von Neumann subalgebras, while in the case of an acting standard $\lambda$-lattice $\Cal G$  (or of an acting rigid $C^*$-tensor category)  
one generally has to deal with {\it quasi-normalizers} (see [BHV15] for related results).

\heading 6. Vanishing cohomology and Connes Embedding  conjecture  \endheading 

In this section, we'll show that Connes Approximate Embedding (CAE) conjecture for factors of the form $R \rtimes \Gamma$ 
can be reformulated as a vanishing cohomology problem for a certain cocycle action of $\Gamma$. Thus, let 
$\omega$ be an (arbitrary) non-principal ultrafilter on $\Bbb N$ and denote by $R^\omega$ the corresponding 
$\omega$-ultrapower II$_1$ factor, with $R\subset R^\omega$ viewed as constant sequences.

\proclaim{6.1. Proposition} $1^\circ$ $R_\omega:=R'\cap R^\omega$ is a $\text{\rm II}_1$ factor 
whose centralizer in $R^\omega$ is equal to $R$, i.e., $R_\omega'\cap R^\omega=R$. 

\vskip.05in

$2^\circ$ Given any $\theta\in \text{\rm Aut}(R)$ there exists a unitary element $U_\theta\in \Cal N_{R^\omega}(R)$ 
such that $\text{\rm Ad}(U_\theta)_{|R}=\theta$. If $U'_\theta \in \Cal N_{R^\omega}(R)$ is another unitary satisfying 
$\text{\rm Ad}(U'_\theta)_{|R}=\theta$, then $U'_\theta = v U_\theta = U_\theta v'$ for some $v, v' \in \Cal U(R_\omega)$. 
Moreover, if $U\in \Cal N_{R^\omega}(R)$ and one denotes $\theta=\text{\rm Ad}(U)_{|R}\in \text{\rm Aut}(R)$, then $U \in  \Cal U(R_\omega)U_\theta$.

\vskip.05in

$3^\circ$  If $\theta, U_\theta$ are as in $2^\circ$ above, then 
$\text{\rm Ad}(U_\theta)_{|R_\omega}$ implements an  element $\theta_\omega \in \text{\rm Out}(R_\omega)$ 
and an element $\tilde{\theta}_\omega=\text{\rm Ad}(U_\theta)_{|R \vee R_\omega} \in  \text{\rm Out}(R \vee R_\omega)$, with   
$\theta\in \text{\rm Aut}(R)$ outer iff $\theta_\omega$ outer and 
iff $\tilde{\theta}_\omega$ outer. 

\vskip.05in

$4^\circ$ The application $\text{\rm Out}(R) \ni \theta \mapsto \tilde{\theta}_\omega \in \text{\rm Out}(R \vee R_\omega)$ 
is a $1$ to $1$ group morphism whose image has trivial scalar $3$-cocycle, with corresponding cocycle crossed product 
$\text{\rm II}_1$ factor $(R\vee R_\omega) \rtimes \text{\rm Out}(R)$ 
equal to the von Neumann algebra generated in $R^\omega$ by $R\vee R_\omega$ and $\{U_\theta \mid \theta 
\in \text{\rm Aut}(R)\}$ $($thus equal to $\Cal N_{R^\omega}(R)''$ as well$)$. 

\vskip.05in

$5^\circ$ Any free action $\Gamma \curvearrowright^\sigma R$ gives rise to a free cocycle action $\tilde{\sigma}_\omega$ 
of $\Gamma$ on $R\vee R_\omega$, by $\tilde{\sigma}_\omega(g)=\text{\rm Ad}(U_{\sigma(g)})_{|R \vee R_\omega}$, $g\in \Gamma$, 
with corresponding $2$-cocycle $v_\omega^\sigma: \Gamma \times \Gamma \rightarrow \Cal U(R_\omega)$. 

\endproclaim
\noindent
{\it Proof}. $1^\circ$ This is a particular case of (Theorem 2.1 in [P13]). 

$2^\circ$ This is well known (see e.g., [C74]) and is due to the fact that any automorphism of $R$ is approximately inner. 

$3^\circ$ Since $U_\theta$ normalizes $R$, it also normalizes its relative commutant $R'\cap R^\omega=R_\omega$, 
and therefore $R\vee R_\omega$ as well. If the automorphism $\theta_\omega$ it implements on $R_\omega$ is 
inner, say implemented by some $v\in \Cal U(R_\omega)$, then $v^*U_\theta \in R_\omega'\cap R^\omega = R$, 
implying that $\text{\rm Ad}(U_\theta)$ is inner on $R$, i.e., $\theta$ is inner. Similarly, if $\text{\rm Ad}(U_\theta)$  
is inner on $R$, then it is inner on $R_\omega$. Since $R\vee R_\omega\simeq R\overline{\otimes} R_\omega$ 
with $\text{\rm Ad}(U_\theta)$ splitting as a tensor product of its restrictions to $R, R_\omega$, 
one also has that this automorphism is inner iff both restrictions are inner. 

$4^\circ$ The II$_1$ factor $R\vee R_\omega$ has trivial relative commutant in $R^\omega$ 
and so if we denote by $\Cal N$ the unitaries in its normalizer that leave $R$ (and thus also $R_\omega$) 
invariant, then $\Cal G=\Cal N/\Cal U(R)\Cal U(R_\omega)$ is a discrete group 
implementing a cocycle action on $R\vee R_\omega$, 
with $(R \vee R_\omega) \vee \Cal N \simeq (R\vee R_\omega) \rtimes \Cal G$. Also, from the construction of the map 
of $\text{\rm Aut}(R) \ni \theta \mapsto U_\theta \in \Cal N$ and part $3^\circ$, we see that  
this map implements an isomorphism $\text{\rm Out}(R) \simeq \Cal G$. 

$5^\circ$ This part is trivial from $3^\circ$ above. 
\hfill $\square$

\vskip.05in
\noindent
{\bf 6.2. Definition}. A II$_1$ factor $M$ (respectively a group $\Gamma$) 
has the {\it CAE property} if it can be embedded into $R^\omega$ (respectively into the unitary group 
of $R^\omega$).  Note that by a result in [R06], $\Gamma$ has a faithful representation into $\Cal U(R^\omega)$ 
iff $R^\omega$ contains a copy of the left regular representation of $\Gamma$, equivalently $L(\Gamma)  \hookrightarrow R^\omega$.  
Thus,  $\Gamma$ has the CAE property iff $L(\Gamma)$ has the CAE property.

\proclaim{6.3. Theorem} Let $\Gamma \curvearrowright^\sigma R$ be a free action of a countable group $\Gamma$ on the hyperfinite 
$\text{\rm II}_1$ factor $R$. The $\text{\rm II}_1$ factor $R\rtimes_\sigma \Gamma$ has the CAE property if and only 
if the $\Cal U(R_\omega)$-valued $2$-cocycle $v_\omega^\sigma$ vanishes, i.e., iff there 
exist unitary elements $\{U_g \mid g\in \Gamma \}\subset \Cal N_{R^\omega}(R)$ 
that implement $\sigma$ on $R$ and satisfy 
$U_gU_h=U_{gh}$, $\forall g, h\in \Gamma$. 

\endproclaim
\noindent
{\it Proof}. Let $M=R\rtimes_\sigma \Gamma$ with $\{U_g \mid g\in \Gamma \}$ denoting the canonical unitaries implementing $\sigma$. 
If $M$ is embeddable into $R^\omega$, then by using the fact that any two copies of the hyperfinite II$_1$ factor in $R^\omega$ are 
conjugated by a unitary element in $R^\omega$, it follows that 
we may assume the hyperfinite subfactor $R$ in $M=R\rtimes \Gamma$ coincides 
with the algebra of constant sequences in $R^\omega$, with the action $\sigma$ on it being implemented by 
$\{U_g\}_g \subset M\subset R^\omega$. By Proposition 6.1.5$^\circ$ above, this means the $2$-cocycle $v^\sigma_\omega$ 
vanishes. 

Conversely, if $v^\sigma_\omega$ vanishes, then we clearly have $R\rtimes_\sigma \Gamma \hookrightarrow R^\omega$. 
\hfill $\square$

\proclaim{6.4. Corollary} Let $\Gamma$ be a countable group and $\Gamma \curvearrowright^\sigma R^{\overline{\otimes}\Gamma}$ the 
non-commutative Bernoulli $\Gamma$-action with base $R$. Let $H$ be an ICC amenable group $($such as 
the group $S_\infty$ of finitely supported permutations of $\Bbb N$, or the lamp-lighter group $\Bbb Z/2\Bbb Z \wr \Bbb Z)$. 
Then $H \wr \Gamma$ is a CAE group iff $v^\sigma_\omega$ vanishes. 
\endproclaim

\vskip.05in
\noindent
{\bf 6.5. Remarks.} 1$^\circ$ It has been shown in [HaS16] that if two groups $H, \Gamma$ are sofic, 
then their wreath product $H \wr \Gamma$ is sofic as well, so in particular it is CAE. 
Taking $H$ to be an (arbitrary) amenable ICC group $H$, for  which by Connes Theorem one has $L(H)\simeq R$, 
it follows that 
the crossed product II$_1$ factor $R\rtimes_\sigma \Gamma =L(H\wr \Gamma)$ is CAE, 
where $\Gamma \curvearrowright^\sigma R^{\overline{\otimes}\Gamma}\simeq R$ 
is the non-commutative Bernoulli $\Gamma$-action with base $\simeq R$ as in 6.4. Thus, 
the corresponding cocycle $v^\sigma_\omega$ vanishes. 
Equivalently, one can choose $U_{\sigma(g)}\in \Cal N_{R^\omega}(R)$ 
so that to be a representation of $\Gamma$. One can in fact show that these unitaries can be taken so that to also 
normalize the ultrapower of the Cartan subalgebra,  i.e., 
$\{U_g\}_g \subset \Cal N_{R^\omega}(R)\cap \Cal N_{R^\omega}(D^\omega)$, where $D = D_0^{\overline{\otimes} \Gamma}$, 
$D_0$ being the Cartan subalgebra of the base.  

$2^\circ$ Given any $\Gamma \in \Cal V\Cal C$, the wreath product group $S_\infty \wr \Gamma$ is CAE by Corollary 6.4 above, 
and thus $\Gamma$ is a CAE group. 
However, one already knows this, since we have seen in Section 4 that $\Cal V\Cal C$ is contained in 
W$^*_{leq}(\Bbb F_2)$, and $L(\Bbb F_2) \hookrightarrow R^\omega$. But while the class $\Cal V\Cal C$ has a lot of restrictions on it (cf. Sections 3 and 4 in this paper), 
the class of CAE groups is manifestly huge, in fact it may well be that all groups are CAE. 

$3^\circ$ Part $4^\circ$ of Proposition 6.1 above naturally leads to the following: 

\vskip.05in 
\noindent
{\bf 6.6. Problem.} {\it Does the cocycle action 
$\theta \mapsto \pi(\theta)\overset{def}\to{=} \tilde{\theta}_\omega$ of $\text{\rm Out}(R)$ on $R \vee R_\omega$ 
have vanishing $2$-cohomology}? {\it Can $\pi$ be perturbed by inner automorphisms to a genuine action}?  {\it Is it true that $H^2(\text{\rm Out}(R))=1$} ? 

\vskip.05in

Due to its ``huge size'' and properties (e.g., all torsion free elements are conjugate in $\text{\rm Out}(R)$, by [C74]), 
one should have $H^3(\text{\rm Out}(R))=1$. If this is the case, then  perturbing 
the cocycle action $\pi$ by inner automorphisms to a  genuine action would be equivalent to untwisting its $2$-cocyle. 
If the CAE conjecture turns out to hold true, then Theorem 6.3 implies that the restriction of $\pi$ 
above to any countable subgroup  $\Gamma \subset \text{\rm Out}(R)$ that implements a genuine action on $R$ 
has vanishing cohomology. 

Nevertheless, even if this is the case, the entire cocyle action $\pi$ cannot probably be untwisted 
to a  genuine action. One way to prove this would be to show that its restriction to a certain countable subgroup cannot be untwisted. This amounts to saying  
that some free cocycle action $(\theta, v_\theta)$ of a countable group $\Gamma$ on $R$ cannot be untwisted when viewed as an action on $R \vee R_\omega$. 
So to start with, this means $(\theta, v_\theta)$ cannot be untwisted as a cocycle action on $R$. 

But the only exemples of cocycle actions $(\theta, v_\theta)$ on $R$ that 
we know to be ``un-twistable'' are the ones provided by Theorem 3.2, which are amplifications $(\sigma^t, v^t)$ of the Bernoulli $\Gamma$-actions 
$\Gamma \curvearrowright^\sigma R^{\overline{\otimes}\Gamma}$ (as defined in [P01a]). 
But if $\pi$ can be untwisted on $\sigma(\Gamma)$, then all these cocycle actions can actually be untwisted on $R \vee R_\omega$. Indeed, this is because 
any countable subgroup in the normalizer of $R\vee R_\omega$ in $R^\omega$ has a ``huge'' non-separable type II$_1$ fixed point algebra. 
As noted in [P01a], if a genuine action has II$_1$ fixed point algebra, then all its amplifications can be untwisted. Thus, if $\pi$ can be untwisted 
when restricted to $\sigma(\Gamma)\subset \text{\rm Out}(R)$ then it can also be untwisted when restricted to $\sigma^t(\Gamma)\subset \text{\rm Out}(R)$. 

All this shows that the Problem  6.6. may lead to some interesting logic-related considerations, especially 
when  taken together with the CAE conjecture.

\head  References \endhead

\item{[AP17]} C. Anantharaman, S. Popa: ``An introduction to II$_1$ factors'', \newline www.math.ucla.edu/$\sim$popa/Books/

\item{[Bi97]} D. Bisch: {\it Bimodules, higher relative commutants and the fusion algebra associated to a subfactor}, 
The Fields Institute for Research in Mathematical Sciences Communications Series, Vol. {\bf 13} (1997), 13-63.

\item{[BoHV15]} R. Boutonnet, C. Houdayer, S. Vaes: {\it Strong solidity of free Araki-Woods factors}, 
to appear in Amer. J. Math, arXiv:1512.04820

\item{[BrO08]} N. Brown, N. Ozawa:  ``$C^*$-algebras and finite dimensional-approximations'', Amer. Math. Soc. 
Grad. Studies in Math. {\bf 88}, 2008. 

\item{[C74]} A. Connes: {\it Outer conjugacy classes of automorphisms of factors}, 
Ann. Ecole Norm. Sup., {\bf 8} (1975), 383-419.  

\item{[C75]} A. Connes: {\it Classification of injective factors},
Ann. of Math., {\bf 104} (1976), 73-115.

\item{[C80]} A. Connes: {\it Correspondences}, handwritten notes, 1980. 

\item{[CJ84]} A. Connes, V.F.R. Jones: {\it Property T 
for von Neumann algebras}, Bull. London Math. Soc. {\bf 17} (1985), 
57-62. 

\item{[CT76]} A. Connes, M. Takesaki: {\it The flow of weights of factors of type} III, Tohoku Math. Journ. {\bf 29} (1977), 473-575.

\item{[D92]} K. Dykema: {\it Interpolated free group factors}, Pac. J. Math. {\bf 163} (1994), 123-135.

\item{[G04]}  D. Gaboriau: {\it Examples of groups that are measure equivalent to the free group}, 
Ergodic Theory Dynam. Systems {\bf 25} (2005), 1809-1827. 

\item{[GL07]} D. Gaboriau, R. Lyons: {\it A Measurable-Group-Theoretic Solution to von Neumann's Problem}, 
Invent. Math., {\bf 177} (2009), 533-540.

\item{[H79]} U. Haagerup: {\it An example of a non-nuclear} C$^*$-{\it algebra which has the metric approximation property}, Invent. Math. 
{\bf 50} (1979), 279-293.

\item{[HaS16]} B. Hayes, A. Sale: {\it The wreath product of two sofic groups is  sofic}, arXiv:1601.\newline 03286

\item{[Hj04]} G. Hjorth: {\it A lemma for cost attained}, Ann. Pure Applied Logic {bf 143} (2006), 87-102. 

\item{[I06]} A. Ioana: {\it Rigidity results for wreath product} II$_1$ {\it factors},  J. Funct. Anal. {\bf 252} (2007), 763-791. 

\item{[I13]} A. Ioana:  {\it Cartan subalgebras of amalgamated free product} II$_1$ {\it factors}. 
(Appendix by A. Ioana and  S. Vaes),  Ann. Sci. Ec. Norm. Super. {\bf 48} (2015),  71-130.

\item{[IPeP05]} A. Ioana, J. Peterson, S. Popa:
{\it Amalgamated Free Products of w-Rigid Factors and Calculation of their
Symmetry Groups},
Acta Math. {\bf 200} (2008), 85-153. 

\item{[J80]} V. F. R. Jones: {\it Actions of finite groups on the hyperfinite type 
II$_1$ factor}, Mem. Amer. Math. Soc., {\bf 237}, 1980. 

\item{[J81]} V. F. R. Jones: {\it A converse to Ocneanu's theorem}, 
Journal of Operator Theory {\bf10} (1983), 61-64. 

\item{[J83]} V.F.R. Jones: {\it Index for subfactors}, Invent. Math. {\bf 72} (1983), 1-25. 

\item{[J99]} V.F.R. Jones: {\it Planar Algebras}, math.OA/9909027

\item{[Ju06]} K. Jung: {\it Strongly $1$-bounded von Neumann algebras},  Geom. Funct. Anal. {\bf 17} (2007), 1180-1200. 

\item{[KV83]} V. A. Kaimanovich, A. M. Vershik: {\it Random Walks on Discrete Groups: Boundary and Entropy}, 
Annals of Probability, {\bf 11} (1983), 457-490.

\item{[KT02]} Y. Katayama, M. Takesaki: {\it Outer actions of a countable discrete amenable group on an AFD factor}, 
in ``Advances in quantum dynamics'' (South Hadley, MA, 2002), 163-171, Contemp. Math., {\bf 335}, AMS, Providence, RI, 2003.

\item{[NT59]} M. Nakamura, Z. Takeda: 
{\it On the extension of finite 
factors} I, II, Proc. Japan Acad., {\bf44} (1959), 149-154, 215-220. 

\item{[NeTu13]} S. Neshveyev, L. Tuset:  ``Compact Quantum Groups and their Representation Categories'',  Cours
Sp\'ecialis\'e, Vol {\bf 20}, Soci\'et\' e Math. de France, Paris (2013). 

\item{[Oc85]} A. Ocneanu: ``Actions of discrete amenable groups on 
factors'', Springer Lecture Notes No. 1138,
Berlin-Heidelberg-New York 1985.

\item{[O03]} N. Ozawa: {\it Solid von Neumann algebras}, Acta Math. {\bf 192} (2004), 
111-117. 

\item{[OP07]} N. Ozawa, S. Popa: {\it On a class of} II$_1$ {\it
factors with at most one Cartan subalgebra}, Annals of Mathematics {\bf 172} (2010),
101-137 (math.OA/0706.3623)

\item{[Pe07]} J. Peterson: L$^2$-{\it rigidity in von Neumann algebras}, Invent. math. {\bf 175} (2009), 417-433. 

\item{[PeT10]} J. Peterson, A. Thom: {\it Group cocycles and the ring of affiliated operators}, Invent. Math. {\bf 185} (2011), 561-592. 

\item{[PiP84]} M. Pimsner, S. Popa: {\it Entropy and index for subfactors}, 
Ann. Sci. Ec. Norm. Super. {\bf 19} (1986),  57-106.

\item{[PiP88]} M. Pimsner, S. Popa: {\it Iterating the basic
construction}, Trans. AMS, {\bf 310} (1988), 127-134. 

\item{[P81]} S. Popa: {\it Maximal injective subalgebras in factors
associated with free groups}, Advances in Math., {\bf 50} (1983),
27-48.

\item{[P83]} S. Popa, {\it Hyperfinite subalgebras normalized by a
given automorphism and related problems}, in ``Proceedings of the
Conference in Op. Alg. and Erg. Theory'' Busteni 1983, Lect. Notes
in Math., Springer-Verlag, {\bf 1132}, 1984, pp 421-433.

\item{[P86]} S. Popa: {\it Correspondences}, INCREST Preprint 56/1986, www.math.ucla.edu/ \newline $\sim$popa/preprints.html

\item{[P89]} S. Popa: {\it Sous-facteurs, 
actions des groupes et cohomologie},
C.R. Acad. Sci. Paris {\bf309} (1989), 771-776. 

\item{[P91]} S. Popa: {\it Classification of amenable subfactors of
type II}, Acta Mathematica, {\bf 172} (1994), 163-255.

\item{[P93]} S. Popa: {\it Approximate innerness and central
freeness for subfactors: A classification result}, in
``Subfactors'', (Proc. Tanegouchi Symposium in
Operator Algebras), Araki-Kawahigashi-Kosaki Editors, World
Scientific 1994, pp 274-293.

\item{[P94a]} S. Popa: {\it An axiomatization of the lattice of higher relative commutants of a subfactor}, 
Invent. Math.,  {\bf 120} (1995), 427-445.

\item{[P94b]}  S. Popa, {\it Symmetric enveloping algebras,
amenability and AFD properties for subfactors}, Math. Res.
Letters, {\bf 1} (1994), 409-425.

\item{[P97]} S. Popa: {\it Some properties of the symmetric enveloping algebras
with applications to amenability and property T}, Documenta Mathematica, {\bf 4} (1999), 665-744.

\item{[P01a]} S. Popa: {\it Some rigidity results for
non-commutative Bernoulli shifts}, J. Fnal. Analysis {\bf 230}
(2006), 273-328 (MSRI preprint No. 2001-005). 

\item{[P01b]} S. Popa: {\it On a class of type} II$_1$ {\it factors with
Betti numbers invariants}, Ann. of Math {\bf 163} (2006), 809-899
(math.OA/0209310; MSRI preprint 2001-024).

\item{[P03]} S. Popa: {\it Strong Rigidity of }  II$_1$ {\it Factors
Arising from Malleable Actions of $w$-Rigid Groups} I, Invent. Math.,
{\bf 165} (2006), 369-408. 

\item{[P06a]} S. Popa: {\it On the superrigidity of malleable
actions with spectral gap},  J. Amer. Math. Soc. {\bf 21}
(2008), 981-1000 (math.GR/0608429).

\item{[P06b]} S. Popa: {\it On Ozawa's Property for Free Group Factors},
Int. Math. Res. Notices (2007) Vol. {\bf 2007}, article ID rnm036,
10 pages, doi:10.1093/imrn/rnm036 published on June 22, 2007
(math.OA/0608451)

\item{[P13]} S. Popa: {\it Independence properties in subalgebras of ultraproduct} II$_1$ {\it factors}, Journal of Functional Analysis 
{\bf 266} (2014), 5818-5846 (math.OA/1308.3982)

\item{[PS01]} S. Popa, D. Shlyakhtenko: {\it Universal properties of}
$L(F_\infty)$ {\it in subfactor theory}, Acta Mathematica, {\bf
191} (2003), 225-257.

\item{[PSV15]} S. Popa, D. Shlyakhtenko, S. Vaes:  {\it Cohomology and} $L^2$-{\it Betti numbers for subfactors and quasi-regular inclusions}, 
math.OA/1511.07329 , to appear in the International Mathematical Research Notices. 

\item{[PV12]} S. Popa, S. Vaes: {\it Unique Cartan decomposition for} II$_1$ 
{\it factors arising from arbitrary actions of free groups}, Acta Mathematica, {\bf 194} (2014), 237-284 

\item{[PV14]} S. Popa, S. Vaes: {\it Representation theory for subfactors}, $\lambda$-{\it lattices and C$^*$-tensor categories},
Commun. Math. Phys. {\bf 250} (2015), 1239-1280. 

\item{[R91]} F. Radulescu: {\it The weak closure of the group algebras associated to free groups are stably isomorphic},  
Comm. Math. Physics {\bf 156} (1993), 17-36.

\item{[R06]} F. Radulescu: {\it The von Neumann algebra of the non-residually finite Baumslag group $\langle a,b \mid ab^3a^{-1}=b^2\rangle$ 
embeds into $R^\omega$}. Hot topics in operator theory, 173-185, Theta Ser. Adv. Math., {\bf 9}, Theta, Bucharest, 2008. 

\item{[Su80]} C.E. Sutherland: {\it Cohomology and extensions of von Neumann algebras }, Publ. RIMS, 
{\bf 16} (1980), 135-176. 

\item{[T79]} M. Takesaki: ``Theory of operator algebras'' I., Springer-Verlag, New York-Heidelberg, 1979.

\enddocument